\documentclass[10pt]{article}
\usepackage{graphicx} 

\usepackage{cite,amssymb,graphicx,caption,lipsum,amsmath}
\usepackage{amsfonts}
\usepackage{enumerate}
\usepackage{epsfig}
\usepackage{xcolor}
\usepackage{psfrag}
\usepackage{hyperref}
\usepackage{babel}
\usepackage{lineno}
\usepackage{chemarr}
\usepackage{rotating}
\usepackage{graphicx}
\usepackage{tcolorbox}

\allowdisplaybreaks

\usepackage[margin=1in]{geometry}

\usepackage{rotating}

\usepackage[colorinlistoftodos]{todonotes}
\usepackage{pgfplots}
\usepackage{stackengine}
\usepgfplotslibrary{fillbetween}
\pgfplotsset{compat=newest}
\usetikzlibrary{decorations.markings, patterns}

\usepackage{tikz}
\usetikzlibrary{math}
\usepackage{pgfplots}
\usepackage{stackengine}
\usepgfplotslibrary{fillbetween}
\pgfplotsset{compat=newest}
\usetikzlibrary{decorations.markings, patterns}

\newtheorem{lemma}{Lemma}
\newtheorem{prop}{Proposition}
\newtheorem{remark}{Remark}

\definecolor{Bluish}{rgb}{0.,0.,0.5}
\definecolor{Reddish}{rgb}{0.5,0.,0.}
\definecolor{PiVD_color}{RGB}{48,136,163}
\definecolor{PiCdD_color}{RGB}{112,188,123}
\definecolor{PiCCd_color}{RGB}{9,13,55}

\def\eps{\varepsilon}
\def\rme{\mathrm{e}}

\def\ds{\displaystyle}

\usepackage[colorinlistoftodos]{todonotes}
\usepackage[normalem]{ulem}

\definecolor{ColorTomasL}{rgb}{0.9,0.2,0.1}

\definecolor{ColorJuanC}{rgb}{0.1,0.3,0.9}

\definecolor{ColorSanti}{rgb}{0.9,0.,0.4}

\definecolor{ColorJosep}{rgb}{1,0.1,0.1}

\usepackage{graphicx}
\graphicspath{{./figures/}}

\usepackage{enumerate}
\usepackage{color}
\usepackage{comment}
\usepackage{hyperref}

\newcommand{\R}{\ensuremath{\mathbb{R}}}

\usepackage{tikz}

\newcommand{\CV}{C_{\!_V}}
\newcommand{\CD}{C_{\!_D}}
\newcommand{\CDV}{C_{\!_{D\!V}}}
\newcommand{\tDV}{t_{\!_{D\!V}}}
\newcommand{\tC}{t_{\!_C}}
\newcommand{\tD}{t_{\!_D}}
\newcommand{\tV}{t_{\!_V}}

\newcommand{\IV}{I_{\!_V}}
\newcommand{\IDV}{I_{\!_{D\!V}}}

\newcommand{\cbeta}{\frac{\partial C}{\partial \beta}}
\newcommand{\cvbeta}{\frac{\partial \CV}{\partial \beta}}
\newcommand{\cdbeta}{\frac{\partial \CD }{\partial \beta}}
\newcommand{\cdvbeta}{\frac{\partial \CDV }{\partial \beta}}
\newcommand{\vbeta}{\frac{\partial V}{\partial \beta}}
\newcommand{\Dbeta}{\frac{\partial D}{\partial \beta}}

\newcommand{\vfinal}{V_{\textit{f}}}
\newcommand{\Dfinal}{D_{\textit{f}}}

\newcommand{\infrate}{\iota}

\newcommand{\PiCDD}{\Pi_{_{\CD\!D}}}
\newcommand{\PiCCD}{\Pi_{_{C \CD}}}
\newcommand{\PiVD}{\Pi_{_{V\!D}}}

\newcommand{\etaD}{\beta \eta}

\usepackage{authblk}

\title{Quasineutral multistability in an epidemiological-like model for defective-helper betacoronavirus infection in cell cultures}

\date{\vspace{-5ex}}
\author[1,2,*]{Juan C. Mu\~noz-S\'anchez}
\author[3,4,5,6,*]{J. Tom\'as L\'azaro}
\author[1]{Julia Hillung}
\author[1]{Mar\'ia J. Olmo-Uceda}
\author[5,6]{Josep Sardany\'es}
\author[1,6,7,+]{Santiago F. Elena}
\affil[1]{Institute for Integrative Systems Biology (I$^2$SysBio), CSIC-Universitat de Val\`encia, Paterna, 46980 Val\`encia, Spain}
\affil[2]{Departament de F\'isica Te\`orica, Universitat de Val\`encia, Burjassot, 46100 Val\`encia, Spain}
\affil[3]{Departament de Matem\`atiques, Universitat Polit\`ecnica de Catalunya (UPC), 08028 Barcelona, Spain}
\affil[4]{Institute of Mathematics, UPC-BarcelonaTech (IMTech), 08028 Barcelona, Spain}
\affil[5]{Centre de Recerca Matem\`atica (CRM), Cerdanyola del Vall\`es, 08193 Barcelona, Spain.}
\affil[6]{Dynamical Systems and Computational Virology, CSIC Associated Unit CRM-I$^2$SysBio, Spain}
\affil[7]{Santa Fe Institute, Santa Fe, NM 87501, USA}
\affil[*]{Equal contribution}
\affil[+]{Correspondence: santiago.elena@csic.es}

\begin{document}

\maketitle

\begin{abstract}  
It is well known that, during replication, RNA viruses spontaneously generate defective viral genomes (DVGs). DVGs are unable to complete an infectious cycle autonomously and depend on coinfection with a wild-type helper virus (HV) for their replication and/or transmission. The study of the dynamics arising from a HV and its DVGs has been a longstanding question in virology. It has been shown that DVGs can modulate HV replication and, depending on the strength of interference, result in HV extinctions or self-sustained persistent fluctuations. Extensive experimental work has provided mechanistic explanations for DVG generation and compelling evidences of HV-DVGs virus coevolution. Some of these observations have been captured by mathematical models. Here, we develop and investigate an epidemiological-like mathematical model specifically designed to study the dynamics of betacoronavirus in cell culture experiments. The dynamics of the model is governed by several degenerate normally hyperbolic invariant manifolds given by quasineutral planes - \emph{i.e.}, filled by equilibrium points. Three different quasineutral planes have been identified depending on parameters and involving: (\emph{i}) persistence of HV and DVGs; (\emph{ii}) persistence of non-infected cells and DVG-infected cells; and (\emph{iii}) persistence of DVG-infected cells and DVGs. Key parameters involved in these scenarios are the maximum burst size ($B$), the fraction of DVGs produced during HV replication ($\beta$), and the replication advantage of DVGs ($\delta)$. More precisely, in the case $0<B<1+\beta$ the system displays tristability, where all three scenarios are present. In the case $1+\beta < B < 1+\beta + \delta$ this tristability  persists but attracting scenario (\emph{ii}) is reduced to a well-defined half-plane. For $B > 1 + \beta + \delta$, the scenario (\emph{i}) becomes globally attractor. Scenarios (\emph{ii}) and (\emph{iii}) are compatible with the so-called self-curing since the HV is removed from the population. Sensitivity analyses indicate that model dynamics largely depend on DVGs production rate ($\beta$) and their replicative advantage ($\delta$), and on both the infection rates and virus-induced cell deaths. Finally, the model has been fitted to single-passage experimental data using an artificial intelligence methodology based on genetic algorithms and key virological parameters have been estimated.
\end{abstract}

\today

\section*{Highlights}
\begin{itemize}
    \item We provide a mathematical model for helper-defective virus dynamics in cell cultures.
    \item Dynamics show multistability governed by quasineutral manifolds.
    \item The production rate of DVGs, virus infection rates and virus-induced cell death largely influence dynamics.
    \item The model successfully fits data from single infection, virus accumulation experiments.
\end{itemize}

\newpage 
\tableofcontents

\section{Introduction}
RNA viruses can quickly adapt and trigger epidemics by crossing species barriers due to their high mutation rates, fast replication, and large population sizes \cite{Duffy2008,Sanjuan2010,Belshaw2011}. However, a high mutation rate is a double-edge sword, as many mutations during infection result in defective viral genomes (DVGs). DVGs cannot complete the infectious cycle by themselves thus depending on the viral proteins synthesized by a wild-type helper virus (HV). The generic term DVG include point mutations, hypermutated genomes, deletions, insertions, and genomic reorganizations \cite{Vignuzzi2019}. Huang and Baltimore coined the term defective interfering particles (DIPs) for viral particles containing DVGs and normal structural proteins encoded by the HV \cite{Huang1970}. DIPs were first identified in the late 1940s by Von Magnus and Gard~\cite{vonMagnus1947} based on the negative impact they exerted on virus accumulation. DIPs rely on a HV for replication, disrupting HV accumulation and impacting viral pathogenesis~\cite{Xu2017,Genoyer2018}. Several studies have shown that viruses rich in DIPs reduce virulence \cite{Cave1985}, induce high interferon levels \cite{Fuller1980}, aid viral persistence \cite{De1980,Kennedy1982}, and modulate infection capacities as shown for different SARS-CoV-2 strains \cite{Campos2022}.

Recent high-throughput sequencing techniques have uncovered the emergence of a plethora of DVGs within a single infected host \cite{Gribble2021,Hillung2024a,Rangel2023,Jaworski2017}. Notably, these studies have shown that distinct subsets of prevalent DVGs recurrently appear along time, suggesting intricate dynamics within the viral population. These dynamics encompass competition, and possibly compensation or cooperation, among various DVGs. Positive selection favours the most competitive DVG variants, indicating their relative fitness concerning the HV and other DVGs. Despite the generation of hundreds or even thousands of distinct DVGs during infections, the majority are lost due, among other factors, to population bottlenecks occurring during \emph{in vivo} transmissions among individuals \cite{McCrone2018} or during \emph{in vitro} diluted serial passages \cite{Novella1996}. However, DVGs might persist for long periods of time in immunosuppressed hosts or in those with comorbidities \cite{Zhou2023}, or if a high ratio between viral particles and susceptible cells (a parameter known as the multiplicity of infection, MOI) is experimentally imposed \cite{Stamfer1971,Huang1973,Holland1976}. During infections of \emph{in vitro} cell cultures, as those that have motivated this modelling work, DVGs accumulate when viral populations are repeatedly passed at a high MOI, while they do not accumulate if MOI is low or strong bottlenecks are imposed at each transmission event \cite{Stamfer1971,Huang1973,Felt2022}. The higher the MOI, the more likely DVGs and HV would coinfect the same cell and thus persist in the population \cite{Stamfer1971,Huang1973,Felt2022}. While the molecular recombination mechanisms by which DVGs are generated are well understood, their role in modulating the outcome of viral infections is sometimes unclear. Understanding the characteristics and functions of DVGs is thus crucial for comprehending the complexity of viral infections and for developing strategies to control or mitigate their impact.

DIPs act as true hyperparasites interfering with the HV replication \cite{Holland1976,Damayanti1999}, competing for resources and reducing its accumulation and transmission efficiency \cite{Huang1970,Huang1973,Barret1984,Roux1991,Smith2016}. Experimental \emph{in vitro} studies with cell cultures have shown that DIPs may engage in an arms race with the HV \cite{Horodyski1983,DePolo1986,DePolo1987,Zwart2013}. Shorter genomes earn an advantage in terms of replication speed compared to the HV \cite{Garcia2004}. Additionally, there is evidence of a stronger form of interference in which the DIPs compete more effectively for the viral replication machinery \cite{DePolo1986,Pattanaik1991}. Since their discovery, and given their effect on the accumulation of the HV, DIPs have attracted the attention of researchers as potential antiviral candidates \cite{Marriott2010,Notton2014,Smith2016}, known as therapeutic interfering particles (TIPs). Following this idea, Xiao \emph{et al.} created a TIP by deleting the capsid-coding region of poliovirus~\cite{Xiao2021}. Remarkably, the administration of this TIP to mice triggered a broad antiviral response against diverse respiratory viruses, including enteroviruses, influenza A virus, and SARS-CoV-2. The broad-spectrum antiviral effects of this synthetic TIP was attributed to local and systemic type I interferon responses. Notably, a single dose not only safeguarded animals from SARS-CoV-2 infection but also stimulated the production of SARS-CoV-2 neutralizing antibodies, protecting against reinfection \cite{Xiao2021}. Another successful application of TIPs was reported by Chaturvedi \emph{et al.} \cite{Chaturvedi2021}. Following a synthetic biology approach, these authors generated artificial particles that had significant anti SARS-CoV-2 effect in cell cultures and primary lung organoids, reducing viral accumulation 10- to 100-fold. Furthermore, intranasal application of these TIPs in infected hamsters suppressed the virus by 100-fold in the lungs, reduced pro-inflammatory cytokine expression, and prevented severe pulmonary edema \cite{Chaturvedi2022}. Interestingly, the prediction of SARS-CoV-2 inhibition by a single TIP administration was obtained from a within-host mathematical model based on differential equations~\cite{Chaturvedi2022}. 

Mathematical models describing \emph{trans} interactions between viral genomes are found in the literature \cite{Gao2009,Sardanyes2010}. The dynamics of DIPs have been extensively studied as an extreme case of complementation \cite{Bangham1990,Szathmary1992,Szathmary1993,Frank2000,Sardanyes2010,StaufferThompson2010,Chao2017,Liang2023}. For example, the early work by Szathm\'ary \cite{Szathmary1992}, presented structured deme models to provide a description of the coexistence of virus segments considering HV and DIPs, sensitive and resistant viruses together with DIPs, covirus pairs (\emph{i.e.}, virus that exist as two or more separated particles all of which must be present for the complete replication cycle to occur), and virus–covirus systems. A deeper analysis of the model with HV and DIPs was later developed in \cite{Szathmary1993} by considering cell populations infected by particles differing in number. Different dynamics, such as stable fixed points, periodic orbits, and strange chaotic attractors were identified with this model. Later, Kirkwood and Bangham \cite{Kirkwood1994} developed a differential equations model to analyze a system consisting of well-mixed host cells, HV, and DIPs under serial passage dynamics. The model successfully explained various dynamic behaviors observed in cell culture experiments: fluctuations in virus accumulation during successive passages and self-curing involving the simultaneous extinction of the HV and DIPs as shown experimentally \cite{Jacobson1979,StaufferThompson2010}. More recently, passage experiments of baculoviruses in moth larvae have also provided experimental evidence of chaotic dynamics between HV and DIPs containing long genomic deletions \cite{Zwart2013}.

Among coronaviruses, deletions are the most common type of DVGs \cite{Campos2022,Hillung2024a}, formed through recombination due to conserved homology in specific regions and/or RNA structures \cite{Liao2014,Jennings1983, Saira2013,Poirier2016}. Indeed, SARS-CoV-2 deletion DVGs have been pervasively found both in cell cultures and in patients \cite{Zhou2023,Campos2022}. In this regard, asymptomatic patients tended to have lower DVG loads than symptomatic ones, and highly diverse populations of DVGs were observed after long-term COVID-19 in an immunosuppressed patient, suggesting a relationship between interferon responses and DVGs \cite{Zhou2023}. Given the pervasiveness of DVGs in betacoronavirus populations \cite{Zhou2023} and the lack of evidences supporting their sustained interference activity (and hence their potential development as TIPs), Hillung \emph{et al.} \cite{Hillung2024a} performed long-term \emph{in vitro} evolution experiments to explore the role of DVGs in betacoronaviruses dynamics of diversification and evolution. Two virus models were chosen for this study, the human coronavirus OC43 (HCoV-OC43) and the murine hepatitis virus (MHV).

\begin{figure}[h!]
\centering
\includegraphics[width=\textwidth]{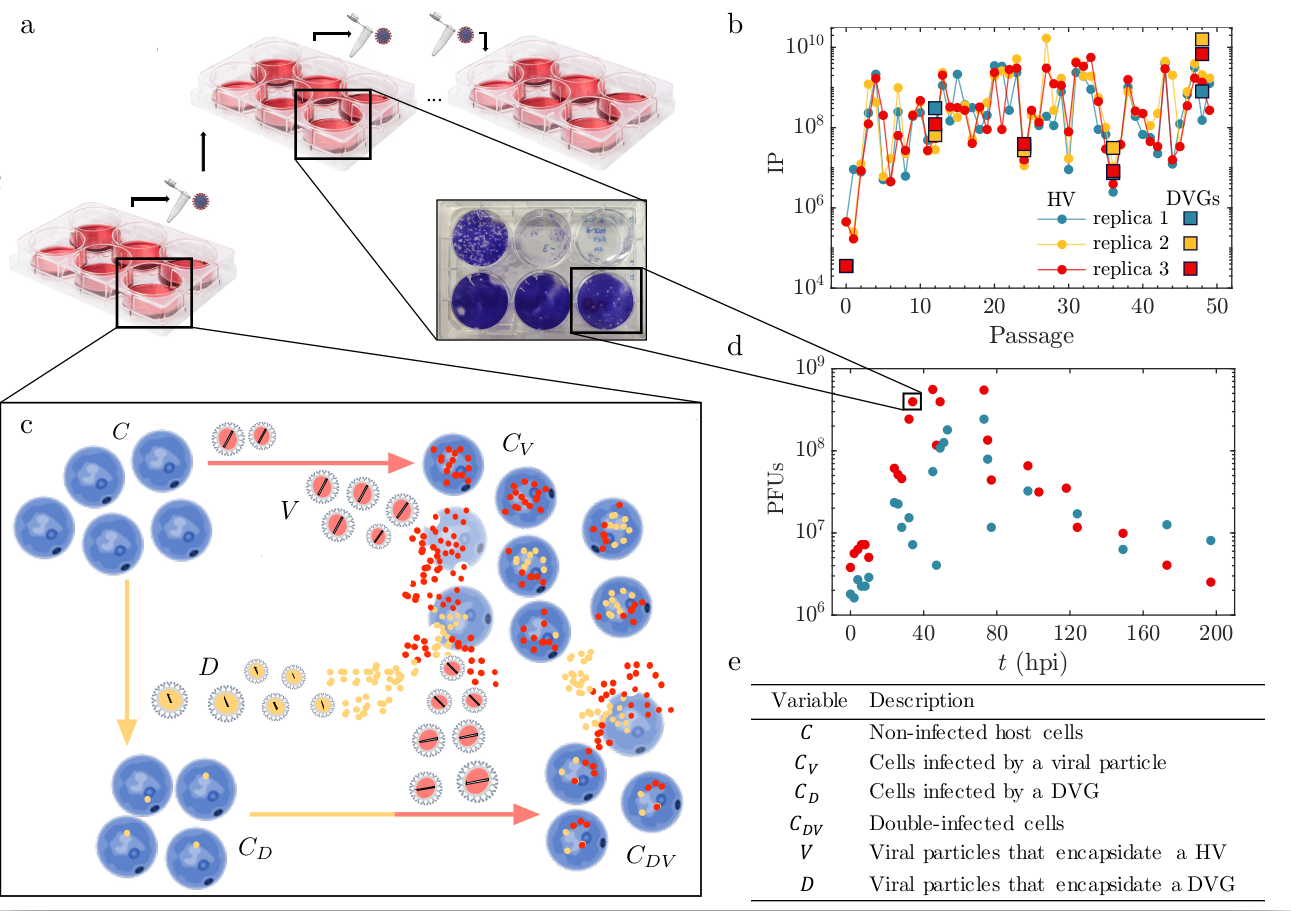}
\caption{(a) Experimental setup for coronaviruses HCoV-OC43 and MHV in cell 6-wells plates~\cite{Hillung2024a}. (b) Passage dynamics for helper virus (HV) and defective viral genomes (DVGs), with large titer fluctuations for the experiments performed in \cite{Hillung2024a}. DVG abundance estimated using \cite{MunozSanchez2024}. (c) Schematic diagram of the dynamical system modelling within- and between-cell virus dynamics for the cell culture  experiments performed in~\cite{Hillung2024b} (see Section~\ref{sec:experiments} and~\cite{Hillung2024a}). The model considers infection of a susceptible population of host cells, $C$, by HV (red) and/or by DVGs (yellow), producing HV only-infected cells $\CV$, DVGs only-infected cells $\CD$, and double-infected cells $\CDV$. Within-cell replication involves the amplification of viral genomes, which produce DVGs at a rate $\beta$. Infected cells are lysed releasing HV and DVGs to the medium. (d) Time series for HCoV-OC43 viral particles infecting BHK-21 cells with different MOIs in cell cultures. (e) State variables of the model.}
\label{fig:diagram}
\end{figure}

To better understand the population dynamics for this experimental system, here we develop and investigate a mathematical model gathering key interactions at the level of within-passage infection dynamics (Fig. \ref{fig:diagram}). The model shows the presence of degenerate normally hyperbolic invariant manifolds (quasineutral invariant manifolds). The existence of such manifolds implies that the orbits in the phase space reach planes of equilibria, depending on parameters. That is, orbits are strongly attracted to a given equilibrium population which is not a single point attractor but a curve or plane filled with equilibria and different initial conditions reach different equilibrium population values. Research on quasineutral planes is rather limited. These neutral surfaces have been characterised in predator-prey models with Holling type III functional responses~\cite{Farkas2005} and in socio-economical models~\cite{Bocso2003}. Moreover, quasineutral states governed by lines or curves of equilibria have been identified in prey-predator models given by differential~\cite{Farkas1984} and partial~\cite{Ferreira2009} differential equations, in Lotka-Volterra competition models \cite{Lin2012}, in strains' competition models of disease dynamics \cite{Kogan2014}, and in models of RNA genomes replication \cite{Sardanyes2018}. More recently, a quasineutral curve was found in a mathematical model for an autocatalytic replicator with an obligate parasite~\cite{Fontich2022}. For this later case, a bistability mechanism determined whether a given initial condition achieved the quasineutral curve or co-extinction. The model investigated in here exhibits a variety of quasineutral objects (mainly planes) displaying, for some parameter combinations, tristability between different scenarios as a function of the initial conditions.

\section{Summary of the experimental results}
\label{sec:experiments}

To investigate the dynamics of defective viral genomes (DVGs) accumulation, we performed serial passages with two betacoronaviruses: HCoV-OC43 in either baby hamster kidney cells (BHK-21) or human large intestine carcinoma cells (HCT-8), and MHV in murine liver cells (CCL-9.1). Passages involved stochastically varying inocula size within two wide but disjoint intervals that can broadly be defined as low and high multiplicity of infection (MOI). At every passage, infectious viral titer (defined as concentration of \emph{plaque forming units}, PFUs/mL) was measured as a proxy to helper virus (HV) accumulation by plaque assays, while the amount and type of the DVGs component of the evolving populations was evaluated at four equidistant passages by high-throughput RNA sequencing (RNA-seq). For illustrative purposes, here we will focus in the case of HCoV-OC43 in BHK-21 cells. In this particular combination, 49 serial passages were performed, and the DVG component evaluated every 12 passages. Fig. \ref{fig:diagram}c shows the dynamics for the first $50$ passages. For the high MOI treatment, the median MOI at the onset of each passage was 25 PFU/cell (IQR: 131.94), while for the low MOI treatment, it was $3.75 \times 10^{-4}$ PFU/cell (IQR: 0.025). See Hillung \emph{et al.} \cite{Hillung2024a} for more experimental details. Each passage involved the replication of the HV and DVGs in a cell culture dish with susceptible host cells. In order to characterize the infection dynamics in the cell culture and test the validity of the model here investigated, a second set of experiments was performed by Hillung \emph{et al.} \cite{Hillung2024b}. Specifically, the dynamics of virus accumulation was monitored within a single infectious passage. To do so, three independent confluent monolayers of BHK-21 cells were inoculated at two viral MOI in the order of units per cell with HCoV-OC43. Then, HV accumulation was evaluated as above at the time points indicated in Fig.~\ref{fig:diagram}b.

\section{Mathematical model}
In this section we introduce the dynamical system modeling experiments by Hillung \emph{et al.} \cite{Hillung2024a,Hillung2024b} summarised in Section \ref{sec:experiments}. The model describes the infection dynamics taking place within a single cell culture dish (Fig. \ref{fig:diagram}a), considering a helper virus (HV), which infects and replicates within a population of susceptible cells, and produces defective viral genomes (DVGs). By assumption, a DVG can only infect a cell but cannot replicate or lysate cells on its own, since it needs the products from the HV. It is also assumed that DVGs have a replication advantage in cells coinfected with the HV. The presence of a DVG inside a cell, by contrast, hinders the replication of the virus, resulting in a fitness penalty. An important assumption made is that viral accumulation dominates against virus' particles decay, \emph{i.e.}, the rate of genomic RNA synthesis is much greater than its degradation rate (previous research on RNA viruses supports this assumption~\cite{Martinez2011}). Furthermore, the model does not consider spontaneous cells' decay and presumes that cell death is driven by virus-induced lysis.  All these assumptions are consistent with the experimental time scale (see below) and aim to focus on the dynamics that arise from the infection process.

The state variables of the model are given by infecting particles, the HV ($V$) and DVGs ($D$), and four different types of cells, including susceptible  ($C$) and infected cells $C_p$, where the subscript indicates the particle (or particles) that has (have) infected the host cell, $p \in \{V, D, DV\}$ (see Fig. \ref{fig:diagram}e). Without affecting any theoretical conclusion and to avoid the use of large numbers (for instance, the number of cells moves around $10^6$) in the propagation of numerical errors, we have scaled all state variables by dividing them by $C(0)$, the initial number of cells. The resulting system is dimensionless and it gives rise to cell variables ranging in the interval $[0,1]$, with $C(0)=1$, and viruses and DVGs taking values in an interval that depends on some parameters of the system.

The infection process is modeled as follows. When a $C$ encounters a $V$ or a $D$ it becomes an infected cell $\CV$ or $\CD$, respectively. This process happens at an infection rate $\infrate$ whose inverse could be considered as the average time between infections. Each type of these infected cells can be coinfected or superinfected with the alternative particle type, resulting in a double-infected cell $\CDV$. $\CD$ cells, as the DVG is not replicating itself, will remain in the same stage until a superinfection with $V$ occurs, if so. The next process will be the replication of HV and DVGs and the lysate of the infected cells. $\CV$ will be lysed resulting into $\eta$ new HV particles and $\beta \eta$ DVGs as a result of errors during the replication process. The total number of particles released from the lysis, or burst size, is $B = \eta (1 + \beta)$. If replication occurs inside a $\CDV$ cell, then the DVG will use the replication machinery of the HV to replicate, competing for common resources and resulting in a decrease of the HV accumulation. The viral production after the $\CDV$ lysis will be reduced to $\eta/\kappa$ with $\kappa > 1$. On the other hand, as DVGs have a replication advantage denoted by $\delta$, if the virus offspring after lysate $\CDV$ is $\eta/\kappa$ then the DVGs produced will be $\eta(\delta + \beta)/\kappa$. This latter term accounts for the HV replication errors resulting into additional DVGs. Thus, the total offspring of the coinfected cells $\CDV$ is $\eta(1 + \delta + \beta )/\kappa$. Under the assumption that the total burst size of $\CV$ and $\CDV$ is the same, as both are the same type of cells, one could derive the penalty coefficient as $\kappa=1+\delta/(1+\beta)$.

For the sake of clarity, the previous processes are represented stoichiometrically with the set of reactions
\begin{eqnarray}
C  + V &\buildrel \infrate \over \longrightarrow& \CV, \\
C + D &\buildrel  \infrate \over \longrightarrow& \CD, \\
\CV + D &\buildrel  \infrate \over \longrightarrow& \CDV, \\
\CD + V &\buildrel  \infrate \over \longrightarrow& \CDV , \\
\CV &\buildrel  \alpha \over \longrightarrow& \eta\, V + \eta\beta\, D , \\
\CDV &\buildrel  \alpha \over \longrightarrow& \frac{\eta}{\kappa}\, V + \frac{\eta(\beta + \delta)}{\kappa} \, D.  
\end{eqnarray}
From these reactions, and using the law of mass action, the following set of autonomous ordinary differential equations (ODEs) is derived:
\begin{eqnarray}
\dot{C} &=& -\infrate \, C (V+D), \label{eq:1}\\[1.2ex]
\dot{\CV} &=& \infrate \, C V - \CV (\infrate \, D+ \alpha), \label{eq:2} \\[1.2ex]
\dot{\CD} &=& \infrate \, (C D - \CD V), \label{eq:3} \\[1.2ex]
\dot{C}_{_{DV}} &=& \infrate \, (\CD V + \CV D) - \alpha \, \CDV, \label{eq:4} \\[1.2ex]
\dot{V} &=& \alpha \, \eta \left( \CV + \dfrac{\CDV}{\kappa} \right) - \infrate \, V (C+\CD), \label{eq:5} \\[1.2ex]
\dot{D} &=& \alpha \, \beta \, \eta \left( \CV + \dfrac{\CDV}{\kappa} \right) + \alpha \, \delta \, \eta \dfrac{\CDV}{\kappa} - \infrate \, D (C + \CV). \label{eq:6}
\label{initial:system}
\end{eqnarray}
Here, $\alpha$, $\beta$, $\delta$, and $\infrate$ are independent parameters, while $\eta$ and $\kappa$ are derived from the stoichiometric relations through
\begin{equation}
\eta=\frac{B}{1 + \beta} \qquad \rm{and} \qquad \kappa=1+\frac{\delta}{1+\beta}.
\label{param:2}
\end{equation}
Notice that the replication advantage of the DVG versus its HV, $\delta$, appears in the term $\alpha \delta \eta \CDV/\kappa$. We assume $\beta \in [0,1]$, $\delta>1$, and $B, \infrate, \alpha > 0$ (see Table \ref{tab:paramdefs}). The assumption $\beta \in [0,1]$ relies on the biological implausibility of $\beta>1$, since DVGs arise due to deletions of the HV genomes. Hence, $\beta$ can be interpreted as the fraction of DVGs produced from HV during cell infection. As we mentioned above, the model assumes that the amplification of both HV and DVGs dominate over their degradation within the experimental time-scales. This assumption is grounded in the experimental measures of accumulation and degradation rates for HCoV-OC43 done in a single passage by Hillung \emph{et al.} \cite{Hillung2024b}. As it can be seen from this system, the scaling performed to the state variables formally affects only the expression of the parameter $\infrate$. These data show that during the first 60 hours post-inoculation (hpi), the effect of degradation was negligible compared to that of production. Only after exhaustion of productive cells, degradation become relevant, though the rate of degradation was still $78.14\%$ slower than the rate of production \cite{Hillung2024b}. In any case, the possible consequences and impact on dynamics of relaxing this assumption are discussed below.

\begin{table}
  \centering
  \caption{Model parameters. All parameters are dimensionless except the infection and lysis rates with inverse time dimensions. HV: helper virus; DVG: defective viral genomes.}
  \begin{tabular}{cll}
  \hline
  Parameter & Description & Range \\
  \hline
  $B$ & Number of particles released after cell lysis (burst-size) & $> 0$ \\
  $\eta$ & Number of HVs produced per cell &$>0$ \\
  $\beta$ & Fraction of DVGs produced per HV due to erroneous replication &$[0,1]$ \\
  $\kappa$ & Replication penalty for HV in cells coinfected with DVGs & $>1$ \\
  $\delta$ & Replication advantage of DVGs & $>1$\\
  $\infrate$ & Infection rate & $>0$\\
  $\alpha$ & Virus infection-induced cell death rate &$>0$ \\
  $m$ & Multiplicity of infection (MOI)& $> 0$\\
  \hline
  \end{tabular}
  \label{tab:paramdefs}
\end{table}

\section{Results}
In the following sections, we compute the equilibria of system Eqs.~\eqref{eq:1}-\eqref{eq:6} and discuss their stability. Phase diagrams for relevant biological parameters are numerically obtained.\footnote{Numerical integrations of the ODEs have been performed using the {Runge-Kutta-Fehlberg-Sim\'o} method of $7^{th}$-$8^{th}$ order with automatic step size control and local relative tolerance $10^{-13}$.} Then, the basins of attraction of equilibria are numerically computed (Section \ref{se:equilibria}). Section \ref{se:general:estimates:vf:Df} provides further results, including those cases with the degradation of viral particles. Section \ref{se:m:larger:1} explores the particular case of multiplicity of infection (MOI) $> 1$ where most cells get infected simultaneously. Next, Section \ref{se:VE:param} contains a sensitivity analysis of the solutions with respect to parameters and to initial conditions. Last but not least, in Section \ref{subsec:fitting} the experimental time series of the HV dynamics in cell cultures generated by Hillung~\emph{et al.}~\cite{Hillung2024b} have been fitted with the mathematical model using artificial intelligence.

\subsection{Planes of equilibria, stability and basins of attraction}
\label{se:equilibria}
As it is commonly done in Dynamical Systems Theory, the study of the dynamics of system~\eqref{eq:1}-\eqref{eq:6} is first based on the computation of its equilibrium solutions and the analysis of their local stability. 

\begin{lemma}[The origin and its local stability]
The origin $\mathcal{O}=(0,0,0,0,0,0)$, the full extinction of all the populations, is always an equilibrium point of the system~\eqref{eq:1}-\eqref{eq:6}. Moreover, for any value of the parameters, the eigenvalues of its jacobian matrix at the origin $\mathcal{O}$ are $0$ (with multiplicity $4$ and semisimple) and $-\alpha$ (double and semisimple as well).
\end{lemma}
Besides the origin, other equilibria are found and discussed in the following proposition. 
\begin{prop}[Planes of equilibria]
\label{lem:equilibria}
System~\eqref{eq:1}--\eqref{eq:6} has three planes formed by equilibrium points (\emph{i.e.}, all the points forming the planes are fixed by the dynamics)
that we label as $\PiCCD$, $\PiCDD$, and $\PiVD$. The origin $\mathcal{O}$ trivially belongs to all of these planes of equilibria, but it does not share their biological interpretation. Because of this, when we refer to these planes the origin will not be considered a part of them, but studied aside. These planes involve different biological equilibrium scenarios and are defined as follows:
\begin{itemize}
\item[(a)] Persistence only of non-infected cells and defective viral genomes (DVG)-infected cells:
\[
\PiCCD = \{(C,\CV,\CD ,\CDV ,V,D) = (C,0,\CD ,0,0,0) \ \big| \ C,\CD \geq 0, \ (C,\CD) \ne (0,0)\}
\]
The spectrum of the jacobian matrix at these points inside the planes is
\begin{equation}
\{ 0\ (\textit{double, semisimple}), \ -\alpha, -\infrate C, \varepsilon_{\pm} \}
\label{Pi:C:CD:vaps}
\end{equation}
where
\begin{equation}
\varepsilon_{\pm} = - \frac{\infrate (C+\CD )+\alpha}{2} \pm \frac{1}{2} \sqrt{  4\alpha C \eta \infrate + \left( \infrate (C+\CD)  -\alpha \right)^2 + \frac{4\alpha \CD \eta \infrate}{\kappa}}.
\label{Pi:C:CD:vaps:mu}
\end{equation}
For any value of the parameters, the discriminant of~\eqref{Pi:C:CD:vaps:mu} is always non-negative, and so all the eigenvalues, for any point in $\PiCCD$, are real. 
 This means that $\PiCCD$ is a so-called normally hyperbolic invariant manifold (NHIM). This plane involves no population of HV (either free or inside cells) or double-infected cells, \emph{i.e.}, a situation of self-curing driven by the quick and efficient outcompetition of the HV by the DVGs. Non-infected cells still remain in the system.

\item[(b)] Persistence only of DVGs and DVG-infected cells:
\[
\PiCDD = \{ (C,\CV,\CD ,\CDV ,V,D)=(0,0,\CD ,0,0,D) \ \big| \ \CD , D \geq 0, \ (\CD ,D)\ne (0,0) \}.
\]
The spectrum of the jacobian at any of its points is
\begin{equation}
\{ 0 \ (\textit{double, semisimple}), \ -\infrate D, -(\infrate D+\alpha), \lambda_{\pm} \},
\label{Pi:CD:D:vaps}
\end{equation}
where
\begin{equation}
\lambda_{\pm} = -\frac{\infrate \CD + \alpha}{2} \pm \frac{1}{2} \sqrt{ (\infrate \, \CD - \alpha)^2 + \frac{4\alpha \CD \eta \infrate}{\kappa}}.
\label{Pi:CD:D:vaps:lambda}
\end{equation}
As in the case above, the discriminant is also non-negative for any choice of the parameters. Hence, all the eigenvalues are real and, 
thus, $\PiCDD$ is also a NHIM. This plane also involves self-curing since, asymptotically, no HV are found as DVGs have steadily and efficiently displaced them from the system.

\item[(c)] Persistence only of free HV and DVGs:
\begin{equation}
\PiVD = \{ (C,\CV,\CD ,\CDV ,V,D)=(0,0,0 , 0, V, D) \ \big| \ V, D \geq 0, \ (V,D)\ne (0,0) \},
\label{def:Pi:VD}
\end{equation}
involving no cells and only HVs and DVGs in the medium. The spectrum of the Jacobian at these points is given by
\begin{equation}
\{ 0,0, -\infrate V, -\alpha, -\infrate (V+D), -(\infrate D+\alpha) \}
\label{Pi:VD:vaps}
\end{equation}
The two (semisimple) $0$-eigenvalues come from the fact that, for any $(V,D)$ the point $(0,0,0,0,V,D)$ is an equilibrium point.
All the eigenvalues are non-positive (and semisimple) and so $\PiVD$ is locally attracting for any $(V,D)\ne (0,0)$.  This plane represents the most common outcome in which HV ends up killing all cells and DVGs are unavoidably present as byproducts of HV replication.
\end{itemize}
\end{prop}
It is not difficult to check that there are no other equilibrium solutions out from the origin and these three planes. 

Let us discuss in more detail the local stability of the equilibrium points contained in the planes. Regarding those on $\PiCDD$, notice that their local stability depends on the sign of the eigenvalue $\lambda_+$ (the rest are all negative, except those $0$ coming from the fact that any point on it with arbitrary $\CD,D$ is also an equilibrium).
Hence, the expression
\[
\lambda_+ \geq 0 \Leftrightarrow \infrate \CD + \alpha \leq \sqrt{ (\infrate \CD - \alpha)^2 + \frac{4\alpha \CD \eta \infrate}{\kappa}}
\]
can be squared (no spurious solutions are introduced since the discriminant is always non-negative), and it leads to
\begin{equation*}
(\infrate\CD + \alpha)^2 - (\infrate \CD - \alpha)^2 \leq \frac{4\alpha \CD \eta \infrate}{\kappa} 
\Leftrightarrow 4\infrate\CD \alpha \leq \frac{4\alpha \CD \eta \infrate}{\kappa} \Leftrightarrow \kappa \leq \eta
\Leftrightarrow 1 + \frac{\delta}{1+\beta} \leq \frac{B}{1+\beta}.
\end{equation*}
That is,
\begin{equation}
\lambda_+ \geq 0 \Leftrightarrow B \geq 1+\beta + \delta.
\label{lambda_plus:nonnegative}
\end{equation}
So, if $B>1+\beta+\delta$ all the points on $\PiCDD$ are unstable and if $B<1+\beta+\delta$, then they are all locally attracting.

Concerning the equilibrium points forming $\PiCCD$ we have, like in the previous case, that their stability depends only on the sign of the eigenvalue $\varepsilon_+$. Thus, squaring again, we obtain
\begin{eqnarray*}
\varepsilon_+ \geq 0 
&\Leftrightarrow& 4\alpha C \eta \infrate + \left( \infrate(C+\CD ) - \alpha \right)^2 + \frac{4\alpha \CD \eta \infrate}{\kappa} \geq \left( \infrate (C+\CD ) + \alpha \right)^2 \\
&\Leftrightarrow& 4\alpha C \eta \infrate + \frac{4\alpha \CD \eta \infrate}{\kappa} \geq \left( \infrate(C+\CD )+\alpha \right)^2 - \left( \infrate(C+\CD ) - \alpha \right)^2 \\
&\Leftrightarrow& 4\infrate\alpha (C+\CD ) \leq 4\alpha C \eta \infrate + \frac{4\alpha \CD \eta \infrate}{\kappa} \Leftrightarrow
(\eta -1) C \geq \left( 1 - \frac{\eta}{\kappa} \right) \CD .
\end{eqnarray*}
Since
\[
\eta - 1 = \frac{B}{1+\beta} - 1, \qquad \quad
1 - \frac{\eta}{\kappa} = 1 - \frac{B}{1+\beta+\delta},
\]
it follows that
\begin{equation}
\varepsilon_+ \geq 0 \Leftrightarrow \left( \frac{B}{1+\beta} -1 \right) C \geq \left( 1 - \frac{B}{1+\beta+\delta} \right) \CD .
\label{mu_plus:positive}
\end{equation}

From conditions~\eqref{lambda_plus:nonnegative} and~\eqref{mu_plus:positive} three possible cases arise (summarized in Table \ref{tab:stability}): (\emph{i}) $B > 1+ \beta + \delta$; (\emph{ii}) $1+\beta < B < 1+\beta + \delta$; and (\emph{iii}) $0<B < 1+\beta$. It is worth remarking that $B$ is one of the parameters that can be better inferred experimentally by counting PFUs as described in Section \ref{se:m:larger:1}. The plane of equilibria $\PiVD$ (taking out the origin) is locally attracting in all three cases, so the following discussion will concern only the planes of equilibria $\PiCCD$ and $\PiCDD$.
\begin{itemize}
\item[(\emph{i})] \textbf{Case} {\boldmath $B > 1+ \beta + \delta$.} On one hand, this case involves $\lambda_+>0$, and so all the equilibrium points forming $\PiCDD$ are, simultaneously, unstable. On the other, $B > 1+ \beta + \delta > 1+ \beta$ and so
\[
\frac{B}{1+\beta} -1 >0 \qquad \textrm{and} \qquad
1-\frac{B}{1+\beta+\delta}<0.
\]
This implies condition \eqref{mu_plus:positive} and so all the equilibria forming $\PiCCD$ become unstable.  Biologically, this is the most commonly expected situation: $B$ contains both HV and DVGs; the greater $\beta + \delta$, the more DVGs are produced at the expense of HV production.

\item[(\emph{ii})]\textbf{Case} {\boldmath $1+\beta < B < 1+\beta + \delta$.} From~\eqref{lambda_plus:nonnegative} it follows that $\lambda_+<0$ and, consequently, all the points constituting $\PiCDD$ are locally attracting. Moreover, one has that
\[
\frac{B}{1+\beta} -1 >0 \qquad \textrm{and} \qquad
1-\frac{B}{1+\beta+\delta}>0.
\]
Thus, from expression~\eqref{mu_plus:positive}, it derives that the points $(C,0,\CD ,0,0,0)\in \PiCCD$ are:
\begin{itemize}
\item[(\emph{iia})] locally attracting if they belong to the half-plane
\begin{equation}
\frac{C}{\CD } < \frac{1 - \frac{B}{1+\beta+\delta}}{\frac{B}{1+\beta} -1 },
\label{c:cd:loc_attracting}
\end{equation}

\item[(\emph{iib})] or unstable (precisely, a saddle) if they fall in
\begin{equation}
\frac{C}{\CD } > \frac{1 - \frac{B}{1+\beta+\delta}}{\frac{B}{1+\beta} -1 }.
\label{c:cd:unstable}
\end{equation}
\end{itemize}

\noindent Notice that condition \eqref{c:cd:loc_attracting} can be equivalently written as
\[
C < \frac{1-\frac{\eta}{\kappa}}{\eta -1} \CD,
\]
with $1<\eta<\kappa$. Since $C,\CD$ are positive, any choice $(\eta_0,\kappa_0)$ inside the open sector bounded by the lines $\eta=1$ and $\eta=\kappa$ (with $\eta,\kappa>1$) gives rise to a region of $(C,\CD)\in (0,+\infty)\times (0,+\infty)$ of attracting equilibrium points. Such regions vary from the void situation (for $\eta=\kappa$), with no point inside, to the one where all the points are attractors (for $\eta \rightarrow 1^+$). In the extreme cases where $\eta = \kappa$, the production of virus per cell would be compensated by the penalty that the HV receives by infecting a $\CD$ cell.

\item[(\emph{iii})] \textbf{Case} {\boldmath $0 < B < 1+\beta$}. Clearly, $B < 1+\beta+\delta$ and therefore all the equilibrium points forming $\PiCDD$ are locally attracting. Besides, one has that $\varepsilon_{+}$ in equation~\eqref{Pi:C:CD:vaps:mu} is strictly negative and, consequently, all the points in $\PiCCD$ are also locally attracting.
\end{itemize}

\begin{table}[h!]
\centering
\caption{Stability analyses: attracting planes in terms of the values of the parameters $B$, $\beta$, and $\delta$.}
\begin{tabular}{cl} \hline\\[-2ex]
Case &  Locally attracting equilibrium points\\ \hline \\[0.1ex]
$B > 1+\beta+\delta$ & $\PiVD$ \\
$1+\beta < B < 1+\beta+\delta$ & $\PiVD$, points in $\PiCCD$ such that~\eqref{c:cd:loc_attracting} holds and $\PiCDD$ \\ 
$0<B<1+\beta$ & $\PiVD, \ \PiCCD$  and $\PiCDD$ \\[1.5ex] \hline
\end{tabular}
\label{tab:stability}
\end{table}

\begin{remark} 
From Proposition \ref{lem:equilibria} and its local stability analysis above, the following statements should be highlighted:
\begin{enumerate}
\item Recall that the planes $\PiVD$, $\PiCCD$, and $\PiCDD$ are highly degenerate in the sense that they are formed by equilibrium points. This fact implies the existence of a couple of zero eigenvalues of the jacobian matrix at any of the equilibrium points.

\item The equilibrium points $(0,0,0,0,V,D)$ constituting the plane $\PiVD$ are all locally attracting for any value of the parameters and for any value of $V$ and $D$. $\PiVD$ is a NHIM (attracting in this case).

\item The local stability of the equilibrium points $(C,0,\CD,0,0,0) \in \PiCCD$ depends on the sign of its eigenvalue $\eps_{+}=\eps_{+}(C,\CD)$. In particular, as seen in (ii) above, there exists a line on $\PiCCD$, namely
\begin{equation*}
\frac{C}{\CD } = \frac{1 - \frac{B}{1+\beta+\delta}}{\frac{B}{1+\beta} -1 },
\end{equation*}
separating those points that are stable from those that are unstable.

\item The local stability of the equilibrium points $(0,0,\CD,0,0,D)\in \PiCDD$ depends only on the value of $\CD$. Precisely, it is attracting if and only if its eigenvalue $\lambda_+=\lambda_+(\CD)<0 \Leftrightarrow B<1+\beta+\delta$, a condition that relies only on the parameters of the system. Moreover, one has that $\lambda_+(\CD)=\eps_+(0,\CD).$

\item Since all the eigenvalues of all the jacobian matrices around the equilibrium points are real, the bifurcations (in stability) that these points undergo are always of \emph{transcritical} type, that is, given by a change in the sign of a real eigenvalue.
\end{enumerate}
\end{remark}

\begin{figure}
\centering
\includegraphics[]{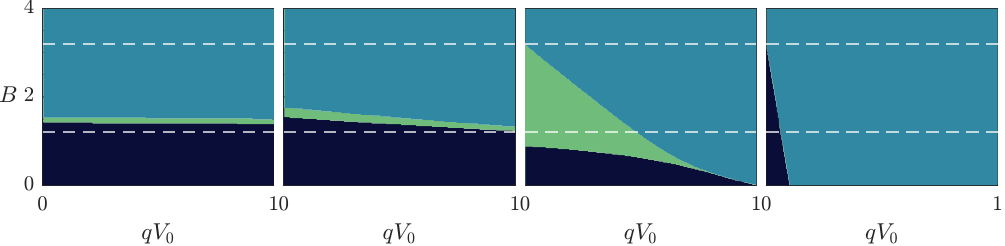}
\caption{Phase diagrams in the plane $(qV_0, B)$. The (approximate) $\omega$-limits are numerically obtained for an orbit with initial conditions $(1,0,0,0,mqV_0, m(1-qV_0))$ in terms of $B$ and the ratio $qV_0 = V(0)/(V(0)+D(0)$. Different colours show the $\omega$-limit given by: \textcolor{PiVD_color}{$\PiVD$} (blue), \textcolor{PiCdD_color}{$\PiCDD$} (olive), and \textcolor{PiCCd_color}{$\PiCCD$} (black). An equilibrium is assumed to be reached when $\|F\|_2 < 10^{-12}$, where $F$ is the vector field of system~\eqref{eq:1}-\eqref{eq:6}. Dashed lines correspond to $B=1+\beta$ and $B=1+\beta+\delta$. The results are displayed for different values of $m$: from left to right $m = 0.01$, $m = 0.1$, $m = 1$, and $m = 10$. In all the panels, we have used $\beta=0.2$, $\delta=2$, and $\infrate/\alpha=10$.}
\label{fig:bifur:Pmax}
\end{figure}

The results shown in Table \ref{tab:stability} lead to some multistability scenarios that can be numerically analysed in terms of the parameters. For example, let us consider the orbits of system~\eqref{eq:1}-\eqref{eq:6} with initial conditions
\begin{equation}
    (C(0),\CV(0),\CD(0),\CDV(0), V(0), D(0)) = (1,0,0,0,mqV_0, m(1-qV_0)),
    \label{eq:condIni_qV0moi}
\end{equation}
and where 
\begin{equation*}
m=\frac{V(0)+D(0)}{C(0)}=V(0)+D(0)
\end{equation*}
is a fixed MOI. Observe that, in this sense, the parameter $qV_0$ provides the virus proportion with this MOI, that is, the ratio $V(0)/(V(0)+D(0))$. 
For any choice $\left(qV_0,  B\right)$, we compute the approximate $\omega$-limit \footnote{Recall that, roughly speaking, the $\omega$-limit of an orbit $\psi(t)$ can be defined as $\omega(\psi)=\lim_{t\rightarrow +\infty} \psi(t)$.} of its corresponding orbit and assign it a colour depending on the plane of equilibria ($\PiVD$, $\PiCDD$, or $\PiCCD$) it achieves. This study leads to different diagrams in terms of the values of $m$. Some of them have been depicted in Fig.~\ref{fig:bifur:Pmax}. Diagrams computed under changes in $\beta$ of around a $10\%-30$\%
provide qualitatively similar results to the ones in Fig.~\ref{fig:bifur:Pmax}. On the contrary, an increase in the value of $\delta$ (from $1.02$ to $5$ in these examples) exhibits substantial growth of the region with final points on the plane $\PiCCD$ (results not shown).

\begin{figure}[h!]
\centering
\includegraphics[width=\textwidth]{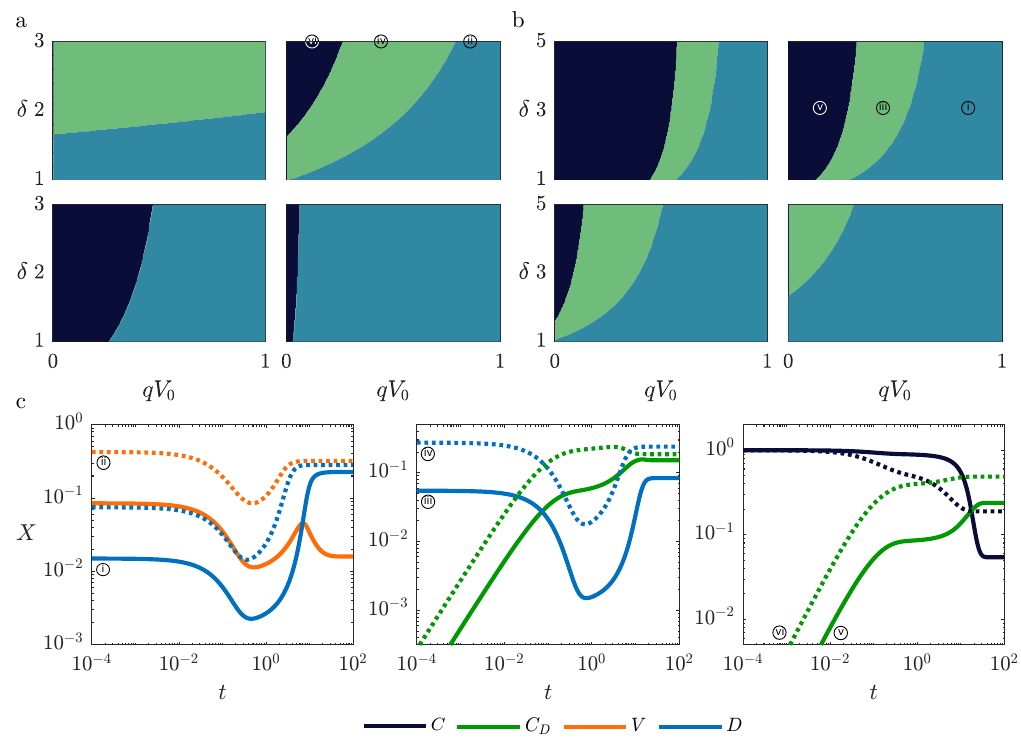}
\caption{ (a) Same as in Fig.~\ref{fig:bifur:Pmax} but fixing $B=1.5$, $\beta=0.2$, $\infrate/\alpha=10$, and tuning $\delta$. (b) Phase diagrams in the plane $(qV_0, \delta)$ for fixed $m = 0.5$ and different values of $B$: from left to right and top to bottom $B = 1.1$, $B=1.5$, $B=2$ and $B=3$. The other parameters are the same as in (a). (c) Examples of trajectories with variables larger than zero at the three planes of equilibria. Numerical values for the different graphs: (\emph{i}) $(qV_0,m) = (0.85, 0.1)$; (\emph{ii}) $(qV_0,m) = (0.85, 0.5)$; (\emph{iii}) $(qV_0,m) = (0.45, 0.1)$; (\emph{iv}) $(qV_0,m) = (0.45, 0.5)$; (\emph{v}) $(qV_0,m) = (0.15, 0.1)$; (\emph{vi}) $(qV_0,m) = (0.15, 0.5)$. The labels (\emph{i}),...,(\emph{vi}) of the time series (at the bottom) correspond to locations in the attraction basin diagrams (top). Other examples of trajectories for all the variables can be found at Fig.~\ref{fig:app:comp:triestability_timeSeries} in Appendix \ref{appendix:timeSeries_examples}.}
\label{fig:bifur:delta}
\end{figure}

Similar explorations can be carried out by setting the values of $\beta$ and $B$ and studying the behaviour when $\delta \geq 1$ varies. As a sample, some of these plots and diagrams are depicted in Fig.~\ref{fig:bifur:delta} for fixed values of the parameters and varying $\delta$.

\begin{figure}[h!]
\centering
\includegraphics[width=\textwidth]{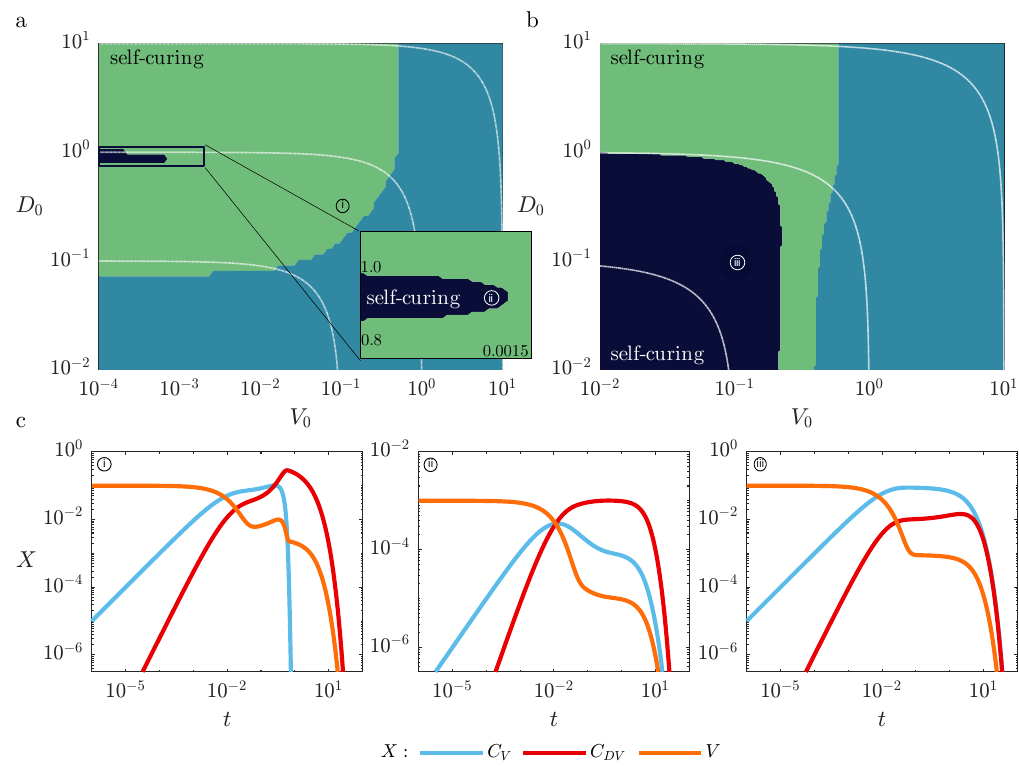}
\caption{Quasineutral multistability for Eqs.~\eqref{eq:1}-\eqref{eq:6} with parameters (a) $B = 10$, $\beta = 0.5$ and $\delta = 20$ and $\infrate/\alpha = 100$ and (b) $B = 1.5$, $\beta = 0.75$ and $\delta = 2$ and $\infrate/\alpha = 100$. The colors show the basins of attraction of each quasineutral plane in the space of initial conditions $(V_0, D_0)$: \textcolor{PiVD_color}{$\PiVD$} (blue), \textcolor{PiCdD_color}{$\PiCDD$} (olive) and \textcolor{PiCCd_color}{$\PiCCD$} (black). The initial conditions for the cell populations are set to zero except for $C_0 = 1$. Lines with constant multiplicity of infection (MOI) are shown in white. (c) Time series showing the self-curing dynamics for three different initial conditions: (\emph{i}) $(V_0,D_0) = (0.10, 0.32)$, (\emph{ii}) $(V_0,D_0) = (0.001, 0.920)$ and (\emph{ii}) $(V_0,D_0) = (0.1, 0.1)$.}
\label{fig:triestability:basins}
\end{figure}

It is  noteworthy to analyse the evolution of the basins of attraction of the planes $\PiCCD$, $\PiVD$, and $\PiCDD$ of Table~\ref{tab:stability} and their biological interpretation in terms of $B$, $1+\beta$, and $1+\beta + \delta$. For instance, the case $B<1+\beta +\delta$ exhibits tristability, in which self-curing and infecting final states simultaneously take place, depending on the initial conditions. The size of the corresponding basins of attraction notably depends on the value of $B$ (with respect to $1+\beta$), the parameters $\beta$ and $\delta$ and the initial infection conditions $(V_0,D_0)$. Recall that cell initial conditions are taken, otherwise indicated, $C(0)=1$ and $\CV(0)=\CD (0)=\CDV (0)=0$. As an example, two different cases are plotted in Fig.~\ref{fig:triestability:basins}, one of them satisfying $1+\beta<B<1+\beta+\delta$ (a) and the other $0<B<1+\beta$ (b), for different initial conditions.  Notice the disparity in the dimensions of the self-curing (both, $\PiCDD$ and $\PiCCD$) basins of attraction. 

\begin{figure}[h!]
\centering
\includegraphics[width=\textwidth]{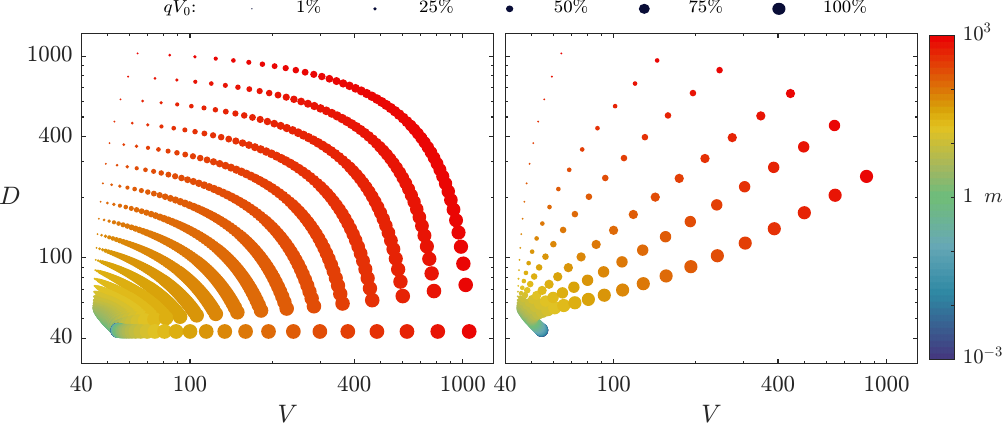}
\caption{ (Left) Equilibria on the plane $\PiVD$ for $B = 100$, $\beta = 0.01$, $\delta = 1.2$, and $\infrate/\alpha = 10$. It can be seen that for $m \gg 1$ both $V$ and $D$ positively correlate with $m$. However, for $m \leq 1 $ this positive correlation with the multiplicity of infection (MOI) disappears as the system becomes more sensitive to parameters and initial conditions. (Right) Selected set of final $(V,D)$-points for six different $qV_0$ values to visually enhance this effect.}
\label{fig:comp_PS}
\end{figure}

Aside from the dynamical interest of the precedent case, the most biologically frequent scenario is certainly the situation where $B > 1+\beta+\delta$, in which a unique globally attractor plane of equilibria $\PiVD$ exists. It is not difficult to prove that any initial infection condition $(V_0,D_0)$ - and cell initial conditions as above - gives rise to a solution which ends, as $t\rightarrow +\infty$, into a final state $(\vfinal,\Dfinal) \in \PiVD$. Clearly, all the cells disappear as the infection-lysis process goes on. In Fig.~\ref{fig:comp_PS}, examples of this final infection state ($\vfinal,\Dfinal)$ for different values of MOI and $qV_0$ are shown.

The existence of these (effective) attracting planes of equilibria in our model is crucial. As already mentioned, this is a consequence of the assumption that virus and cell degradation is much smaller than the amplification of HV and DVGs, and thus can be neglected. Considering degradation of cells and viral particles involves the origin (\emph{i.e.}, total extinction) to be a globally asymptotically stable equilibrium: cells are lysed and degrade, and viral particles will also vanish since no new susceptible cells will remain available for further virus production. After a few hundred hours post-infection (hpi) nothing measurable would remain on the cell culture dishes. Besides, from a mathematical viewpoint, cases with extremely low degradation would be close to the parametric condition involving invariant quasineutral planes. Being extremely close to the structurally unstable condition giving place to the quasineutral objects may involve extremely slow dynamics when the orbits lie close to them (before ending at the origin). Since experimental passages take place during a finite, short time window (between 24 and around 72 hpi) the observed solutions may correspond to a point on some of these three attracting planes or to solutions being trapped in long transients close to these manifolds. As mentioned, one may expect this to occur when parameter values are slightly close to those values giving place to the quasineutral manifolds (see \emph{e.g.}, \cite{Fontich2022}).

\subsection{Cell culture dynamics: estimation of the amount of infective particles}
\label{se:general:estimates:vf:Df}

In this section, we focus our attention on the dynamics leading to asymptotic states, which contain both HV particles and DVGs. Let us assume that $B > 1+ \beta + \delta$ and so the unique attracting plane (filled with equilibrium points) is $\PiVD$ (see Eq.~\eqref{def:Pi:VD}). As already mentioned, we aim to provide estimates for the final quantity of HVs and DVGs obtained after inoculating a new cell culture dish in a given passage. This is the most common experimentally observed case.

Mathematically, the time needed by an orbit to reach an equilibrium \emph{e.g.}, $(0,\ldots,0,V_{\infty}, D_{\infty})$ is infinite, experiments show that at around 72 hpi~\cite{Hillung2024a} cells $C, \CD , \CV$, and $\CDV $ have become extinct and thus one can obtain a good aproximation for the equilibria $(V_{\infty}=\vfinal, D_{\infty}=\Dfinal)$. Strictly speaking, in the model, both $V$ and $D$ would exponentially decrease due to spontaneous degradation. Numerically, inspired by what one can measure in the laboratory, we will consider extinctions if, in absolute value, they are below a (quite small) threshold, say $10^{-10}$ \footnote{Notice that only $C=0$ constitutes a (mathematical) barrier for the values of $C$. It is an invariant manifold and so it cannot be crossed by the flow. This is not the case for the rest of the cell types. This means that at some very particular moments during the numerical integration of the system, a value of $\CD $, $\CV$ or $\CDV $ could be (tiny) negative.} This simplification does not produce relevant numerical consequences (remember that we are computing them for a finite time interval) and accelerates the stabilisation of the values of $V$ and $D$ (as observed experimentally). The orbits we are considering are those having initial conditions of the form
\[
(C(0),\CV(0),\CD(0) ,\CDV(0) ,V(0),D(0))=(1,0,0,0,V_0,D_0),
\]
at time $t=0$ and MOI $m = V_0+D_0$. The integration time is expected to be around 72 hpi, the time invested for a one-passage experiment.

\begin{figure}[h!]
\centering
\includegraphics[width=.86\textwidth]{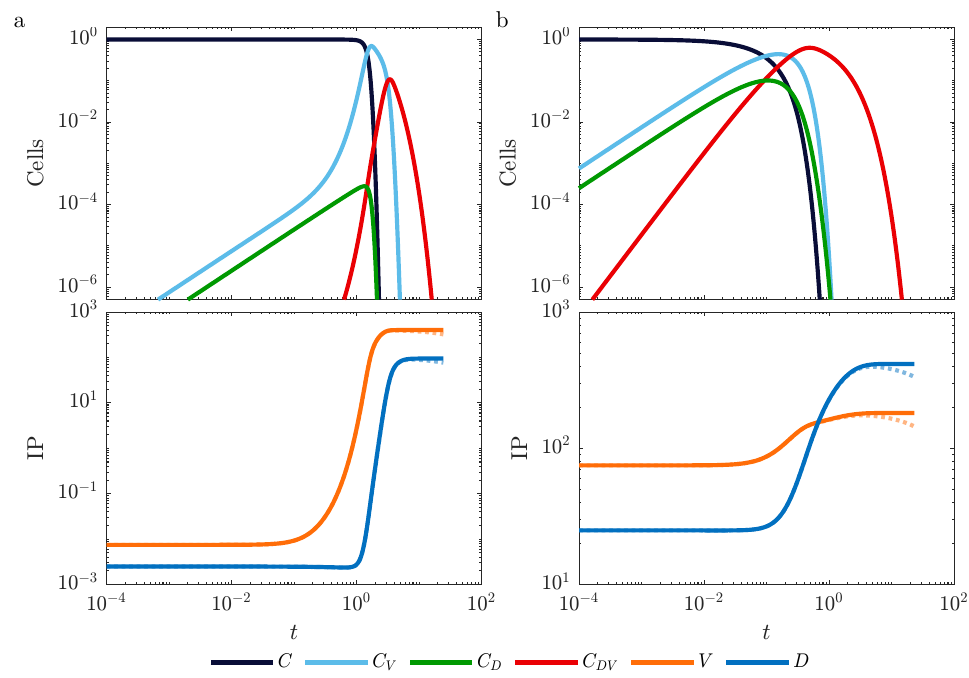}
\caption{Examples of time series for cells and infective particles (IP) at (a) low ($m = 0.01$ and $qV_0 = 0.75$) and (b) high ($m = 100$ and $qV_0 = 0.75$) multiplicities of infection (MOIs). Recall that the quantities are given in units of $10^6$ particles (that is, in such a way that $C(0)$ corresponds to $1$). The extinction times $\tC, \tV$, and $\tDV$ satisfy $0<\tC<\tV<\tDV$. The parameters used are $B = 500$, $\beta = 10^{-6}$, $\delta = 10$ and $\infrate/\alpha = 0.1$. Dashed lines represent the effect of small degradation with $\gamma/\alpha = 0.01$.}
\label{fig:traj_ex}
\end{figure}

Virus infection in a culture dish experiment complies (at any passage) with the following three processes:
\begin{enumerate}
\item \textbf{Infection}. This is the initial state. A certain quantity of viral particles infects the cells. This is noticed in the experimental data (and in the simulated model) by a decrease in the number of PFU/mL observed in the supernatant (outside the cells).

\item \textbf{Cell lysis}. We assume that cells only die as a result of HV infection and replication. The cellular lysis releases HVs and DVG particles into the medium; HV  can be counted from the supernatant by plaque assays (as PFUs). The cell dies afterwards. One could expect the $\textrm{PFU}$-curve to reach a maximum when all the cells have lysed if $m \geq 1$ or to exhibit several waves in the case of $m < 1$, as seen \emph{e.g.}, in \cite{Cuevas2005}.

\item \textbf{Degradation}. Although its effect is not implicitly considered in our model, degradation affects both cells (of any type) and infective particles ($V$ and $D$). For simplicity, one can assume all of them follow a similar exponential law of the form $\rme^{-\gamma t}$, where $\gamma$ is the degradation rate. For example, in the case of HCoV-OC43 one has $\gamma = 0.021 \pm 0.006$ h$^{-1}$ in media containing cellular debris~\cite{Hillung2024b}. From the experimental point of view, it is important to notice that for $m > 1$ the degradation process affects the infective particles in two ways: (\emph{i}) affecting the newly released particles after cell lysis and (\emph{ii}) affecting, from the beginning of the experiment, the $(m - 1)C_0$  particles which do not invade any cell. They constitute a sort of \emph{reservoir}. Since degradation is not taken into account in our model, the final $(\vfinal,\Dfinal)$ should be taken as an overestimation. It is not difficult to see that this difference is of the same order of magnitude as the reservoir's size. For instance, in the case of the HCoV-OC43, only a tiny amount of the initial reservoir would remain infectious after $72$ hpi. This would not change either the order of magnitude of such amounts nor their qualitative behaviour. Examples for low and high MOI time series are shown in Fig.~\ref{fig:traj_ex}. The effect of degradation at a ratio of $\gamma/\alpha = 0.01$ is represented by dashed lines. As said, with a large MOI, the effect of deterioration is greater.

\begin{remark} 
The planes of equilibria $\PiVD$, $\PiCCD$, and $\PiCDD$ are found when decay rates for both viral particles and cells are neglected. However, for very small values of such decay rates the dynamics close to these manifolds is extremely slow and thus the orbits remain close to them for a large amount of time. This involves that, despite the planes no longer exist, the orbits still get trapped there for a long, finite amount of time. See Ref.~\cite{Fontich2022} for an example of this phenomenon for a quasineutral curve.
\end{remark}

The remark above highlights that if the values of the degradation of viral particles and cells remains very low, the orbits may still spend large times close to the regions of the phase space where the planes of equilibria existed. This may involve that such ghost manifolds may be also (temporarily) detected for very small values of decay rates (see \emph{e.g.}, time series with dashed lines in the lower panels of Fig.~\ref{fig:traj_ex}). 
\end{enumerate}

As mentioned above, the aim of this section is to derive insights on the values of $(\vfinal,\Dfinal)$ provided by numerical solutions of the model. Such approximate values become crucial when the inverse fitting problem is tackled: from the quantities of $V$ and $D$ measured at some given times, to estimate the value of some of the parameters of the model providing 'good' fittings of the experimental data. We begin with a more general framework, and afterwards, we perform a second analysis for the particular case of $m > 1$. 

Since in the  numerical simulations, whenever $C, \CV, \CDV$, and $\CD $ are smaller than $10^{-10}$ they are taken as $0$ (extinct), one can define times $\tC, \tV, \tDV$, and $\tD$ such that
\begin{equation}
C(\tC)=0, \qquad \CV(\tV)=0, \qquad \CD (\tD)=0, \qquad \CDV (\tDV)=0.
\label{times_vanishing:cells}
\end{equation}
(see, for instance, Fig.~\ref{fig:traj_ex}). Together with the fact that $\PiVD$ is a globally attracting plane for $B>1+\beta+\delta$, one can expect these times to be ordered as $0< \tC < \tV < \tDV$, with $\tD$ moving around $\tC$ and $\tV$. This has been verified in most of the numerical simulations and it will be assumed hereafter.

To get estimates on $V_f$ and $D_f$, we compare the time evolution of $C$ with the ones of $\CV$ and $\CDV $. We write: 
\begin{eqnarray}
\dot{\CD} &=& \infrate (CD - \CD V) = \infrate CD + \infrate CV - \infrate CV - \infrate \CD V \nonumber \\
&=& \infrate C(V+D) - \infrate V (C+\CD ),\label{cd:c}
\end{eqnarray}
and
\begin{eqnarray}
\dot{\CV} &=& \infrate CV - \infrate \CV D - \alpha \CV = \infrate C V + \infrate C D - \infrate C D - \infrate \CV D - \alpha \CV \nonumber \\
&=& \infrate C(V+D) - \infrate D (C+\CV) - \alpha \CV .\label{cv:c}
\end{eqnarray}
From~\eqref{eq:1} and~\eqref{cd:c} it follows that
\[
\dot{C} + \dot{\CD} = \dot{V} - \alpha \eta \left( \CV + \frac{\CDV }{\kappa} \right)
\]
and, integrating in $[0,t]$, that
\begin{equation}
C(t) - C(0) + \CD (t) - \CD (0) = V(t) - V(0) - \alpha \eta \int_0^t \CV(s) \, ds - \frac{\alpha \eta}{\kappa} \int_0^t \CDV (s) \, ds.
\label{c:cd:1}
\end{equation}
If we introduce the notation
\begin{equation}
\IV = \int_0^{\tV} \alpha \CV(s) \, ds, \qquad \qquad
\IDV = \int_0^{\tDV} \alpha \CDV (s) \, ds,
\end{equation}
evaluate~\eqref{c:cd:1} at $t=\tDV$, take into account that $\tV < \tDV$, and use that $\CD (0)=0$, $\CD(\tDV)=0$, we obtain
\begin{equation}
V(\tDV)=V(0) - C(0) + \eta \IV + \frac{\eta}{\kappa} \IDV.
\label{eq:v:t_dv}
\end{equation}
Notice that for $t>\tDV$ we have that $C(t)=\CD (t)=\CDV (t)=0$ and so $\dot{V}=0$, which makes $V(t)=V(\tDV)$ $\forall t>\tDV$. Therefore, we have $\vfinal=V(\tDV)$, where $\vfinal= \lim_{t\rightarrow +\infty} V(t)$. Analogously, from equations~\eqref{eq:1} and~\eqref{cv:c}, we get
\[
\dot{\CV} + \dot{C} = \dot{D} - \alpha \beta \eta \left( \CV + \frac{\CDV }{\kappa} \right) - \alpha \etaD \frac{\CDV }{\kappa} - \alpha \CV
\]
and
\[
\CV(t) - \CV(0) + C(t) - C(0) = D(t) - D(0) - \alpha(\beta \eta +1) \int_0^t \CV(s) \, ds - \frac{\alpha \eta}{\kappa} (\beta + \delta)\int_0^t \CDV (s) \, ds.
\]
As before, evaluating at $t=\tDV$, having in mind that $\tV < \tDV$ and that $\CV(0)=\CV(\tDV)=0$, $C(\tDV)=0$, it follows that
\begin{equation}
D(\tDV)=D(0) - C(0) + ( \beta \eta +1) \, \IV +
\frac{\eta}{\kappa} (\beta + \delta) \, \IDV.
\label{eq:D:t_dv}
\end{equation}
In a similar manner to the precedent case, we get $\Dfinal=\lim_{t\rightarrow +\infty} D(t)= D(\tDV)$.

Let us define the total number of cells as $C_{\textrm{T}} =C+\CV+\CD +\CDV$. It is straightforward to check that it satisfies the ODE
\begin{equation}
\dot{C_T} = -\alpha(\CV + \CDV)
\end{equation}
and so
\[
-C(0) = -C_{\textrm{T}}(0) + C_{\textrm{T}}(\tDV) = \int_0^{\tDV} \frac{dC}{ds}(s) \, ds = -\alpha \left( \int_0^{\tDV} \CV(s) \, ds + \int_0^{\tDV} \CDV (s) \, ds \right) = -(\IV + \IDV),
\]
from where
\begin{equation}
\IV + \IDV = C(0)
\label{eq:Iv:Idv}
\end{equation}
follows. That is, the initial number of susceptible cells $C(0)$ has been at the end of the process infected either by $V$ or by $V$ and $D$. From this expression, the viral infection-induced death rate $\alpha$ can be estimated as
\[
\alpha = \frac{C(0)}{\ds \int_{0}^{\tV} \CV(t) \, dt + \int_{0}^{\tDV} \CDV(t) \, dt }
\]
or, equivalently,
\[
\alpha = \frac{C(0)}{\ds \int_{0}^{\tDV} (\CV(t) + \CDV(t)) \, dt}.
\]
A large $\alpha$ means that the sum of the two accumulations $\IV$ and $\IDV$ is low, meaning that $\CV$ lyse quickly before defective particles $D$ have time to superinfect them and give place to $\CDV$. On the contrary, $\alpha$ small means that the sum of accumulations is large. One way for this to happen is that cells infected by $V$ (which have already contributed to the accumulation sum $\IV$) are superinfected by a $D$ so that they also contribute to the accumulation $\IDV$.

Following a similar argument, we have that
\[
\dot{\CV} + \dot{\CD} + \dot{C}_{_{CDV}} = \infrate C(V+D) - \alpha (\CV + \CDV).
\]
Hence, integrating in $[0,\tDV]$, and having in mind that
\[
\CV(\tDV)=\CD(\tDV)=\CDV(\tDV)= 0, \qquad \CV(0)=\CD(0)=\CDV(0) =0,
\]
and relation~\eqref{eq:Iv:Idv}, it turns out that the infectivity rate $\infrate$ is given by
\[
\infrate = \dfrac{C(0)}{\ds \int_0^{\tDV} C(t)(V(t)+D(t)) \, dt},
\]
where it gets explicit now that $\infrate$ can be interpreted as the inverse of the mean time among infections.

Combining the expressions above for $\alpha$ and $\infrate$ as
\begin{equation}
    \frac{\infrate}{\alpha} = \dfrac{\ds \int_{0}^{\tDV} (\CV(t) + \CDV(t)) \, dt}{\ds \int_0^{\tDV} C(t)(V(t)+D(t)) \, dt}
    \label{lameua}
\end{equation}
helps to get an interesting insight regarding DVGs accumulation. The greater the ratio $\infrate/\alpha$, the longer $\CV$ cells are available for superinfection by $D$, whereas if HV kills cells quickly (big $\alpha$) and the superinfection is slow (small $\infrate$), $D$ has fewer chances to superinfect and replicate. Therefore, the $\infrate/\alpha$ can be seen as the efficiency of DVG production \emph{i.e.}, efficiency of generating virus-producing infected cells in presence of $D$.

From relation~\eqref{eq:Iv:Idv}, expressions~\eqref{eq:v:t_dv}-\eqref{eq:D:t_dv}, the definitions of $B$, $\kappa$, and denoting $C_0=C(0)$, $V_0=V(0)$, and $D_0=D(0)$, we get
\begin{eqnarray*}
\vfinal + \Dfinal &=& V_0 + D_0 - 2 C_0 + \left( (\beta+1) \eta + 1 \right) \IV + \frac{\eta}{\kappa} (\beta + \delta +1 ) \IDV \\
&=& (m -2) C_0 + (B +1) \IV + B \IDV \\
&=& (m -2) C_0 + B \left( \IV + \IDV \right) + \IV \\
&=& (m - 2 + B) C_0 + \IV.
\end{eqnarray*}
That is,
\begin{equation}
\vfinal + \Dfinal = (m - 2 + B) C_0 + \IV.
\label{eq:vfinalDfinal}
\end{equation}
Equivalently, if in~\eqref{eq:vfinalDfinal} we use equality~\eqref{eq:Iv:Idv}, it follows that
\begin{equation}
\vfinal + \Dfinal = (m - 1 + B) C_0 + \IDV.
\label{eq:vfinalDfinal:2}
\end{equation}
Moreover, from the inequalities
$0 \leq \IV \leq \IV + \IDV = C_0$,
it turns out that
\begin{equation*}
(m - 2 + B)C_0 \leq \vfinal + \Dfinal \leq (m - 1 + B) C_0
\label{eq:ineq:vf:Df}
\end{equation*}
or, equivalently,
\begin{equation}
B + m - 2 \leq \frac{\vfinal+\Dfinal}{C_0} \leq B + m - 1.
\label{eq:ineq:vf:Df:C0}
\end{equation}
Notice that $(\vfinal + \Dfinal)/C_0$ would correspond to the expected next-passage MOI in the case passages are sequentially concatenated, with the same number of initial cells $C_0$ and no dilution applied.

\subsection{Case $m > 1$: all cells are simultaneously infected}
\label{se:m:larger:1}

The expressions obtained in Section~\ref{se:general:estimates:vf:Df} for $\vfinal$ and $\Dfinal$ can be slightly complemented in the case where the MOI is $m > 1$. This assumption makes reasonable to expect a unique wave of cell infections: no HV or DVGs released after the first cell lysis wave will have chances of infecting new cells (of types $C$ and $\CV$), since the latests will all have vanished. Here we present a simplified approach under this assumption, which accounts for the resulting particles. Despite the fact that degradation is not taken into consideration either, the particles' reservoir is explicitly computed.

The results are illustrated with a couple of scenarios, which depend on the relative proportion between initial HV and DVGs. It is not pretended (by far) to be exhaustive and to cover all the possibilities, but to compute approximate expected values for $\vfinal$ and $\Dfinal$. Remind that $V_0, D_0$, and $C_0$ denote the initial conditions $V(0)$, $D(0)$, and $C(0)$, respectively, and that we study the case $C(0)=C_0$ and $\CV(0)=\CD (0)=\CDV (0)=0$. Furthermore, we assume that $B > 1+\beta+\delta$ holds and that, consequently, the $\omega$-limit of such orbits falls into $\PiVD$.

As mentioned above, since $m>1$ not all HV and DVGs are statistically expected to invade a cell. This means that a certain amount of both could (hypothetically) remain outside the cells, without undergoing any infection. This set constitutes a reservoir. We will distinguish between those HVs and DVGs released from infected cells after lysis $(V^{(i)},D^{(i)})$, and those staying in the reservoir $(V^{(r)},D^{(r)})$. Moreover, both types of particles should infect with the same probability (the difference is in their replication speed). Thus, we obtain the following cases:
\begin{enumerate}
\item Case $V_0>C_0$, $D_0=0$:
this case corresponds to the very first step of the experiment. We assume an initial population of HV particles and no DVGs in the cell culture. Since $m>1$ one can expect a one-infection wave. The new particles $V$, $D$ come from two sources:
\\
\noindent (\emph{i}) Infection plus lysis: $V_0=C_0$ and so $\vfinal^{(i)}=\eta C_0$, $\Dfinal^{(i)}= \beta \eta C_0$, the latter one obtained from replication errors of the HV. This process of infection and lysis is assumed to happen almost simultaneously and instantaneously.
\\
\noindent (\emph{ii}) From the reservoir: $\vfinal^{(r)}=V_0 -C_0=(m-1)C_0$ and $\Dfinal^{(r)}=0$. 
Thus,
\[
\vfinal = \vfinal^{(i)} + \vfinal^{(r)} = (\eta + m -1) C_0, \qquad \quad
\Dfinal = \Dfinal^{(i)}= \beta \eta C_0,
\]
and therefore
\begin{equation}
\vfinal + \Dfinal = \left(B + m -1 \right) C_0,
\label{moi:larger:1:vfplusDf:case1}
\end{equation}
where we have used that $B=\eta(1+\beta)$. Notice that in this case we have always $\vfinal \geq \Dfinal$. Indeed, by \emph{reductio ad absurdum},
\begin{equation*}
\Dfinal > \vfinal \Leftrightarrow \beta \eta C_0 > (\eta + m-1) C_0 \Leftrightarrow (\beta-1) \eta > m -1.
\end{equation*}
Since $m > 1$ this implies that $\beta - 1>0 $, which is a contradiction with the fact that $\beta \in [0,1]$.

\item Case $V_0>C_0>D_0>0$. As mentioned, HV and DVGs would enter the cells with the same probability. We detail this time evolution in several steps.

\noindent (\emph{a}) $V$ and $D$ particles infect the $C=C_0$ initial cells proportionally to their density, producing $\CV$ and $\CD$, respectively. Hence, after this first invasion, we get
\[
C=0, \qquad \CV=\frac{V_0}{V_0+D_0} C_0 = \frac{V_0}{m}, \qquad 
\CD=\frac{D_0}{V_0+D_0} C_0 = \frac{D_0}{m},
\]
with
\[
V_a=V_0 \left( 1 - \frac{1}{m} \right) > 0, \qquad D_a = D_0 \left( 1 - \frac{1}{m} \right)
\]
remaining in the medium. Notice that the number of remaining viruses $V_a$ is greater than the number of $D$-infected cells $\CD$. Indeed,
\[
V_a > \CD \Leftrightarrow V_0 \left( 1 - \frac{1}{m} \right) > \frac{D_0}{m} \Leftrightarrow V_0 > C_0.
\]
\noindent (\emph{b}) In the second step, $\CD$ cells become $\CDV$ through the invasion by HV and $\CV$ become $\CDV$ by the one of DVGs. Moreover, a reservoir (constituted only by $V$) remains. That is,
\[
(b_1) \ \ \CD + V_a \stackrel{\infrate}{\longrightarrow} \CDV, \qquad \qquad 
(b_2) \ \ \CV + D_a \stackrel{\infrate}{\longrightarrow} \CDV.
\]
The performance of $(b_1)$ leads to $\CD=0$, $\CDV^{(1)}=\CD=D_0/m$ and a remainder of HVs
\[
V_1^{(r)}= V_a - \CD = V_0 \left( 1 - \frac{1}{m} \right) - \frac{D_0}{m} = V_0 - \frac{V_0+D_0}{m} = V_0 - C_0.
\]
Regarding $(b_2)$, $\CV$-cells become $\CDV$ (by the invasion of defectives) or remain $\CV$. Precisely, it gives rise to
\[
\CDV^{(2)}= D_a = D_0 \left( 1-\frac{1}{m} \right), \qquad \CV^{(2)}=\left( \frac{V_0}{m} - D_0 \left( 1 - \frac{1}{m} \right) \right) = C_0 - D_0,
\]
with no reservoir of DVGs. 

So, at the end of steps $(a)-(b)$, regarding cells, we have $C=0$, $\CD=0$, $\CV=\CV^{(2)}=C_0-D_0>0$, and
\[
\CDV=\CDV^{(1)}+\CDV^{(2)}= \frac{D_0}{m} + D_0 \left( 1 - \frac{1}{m} \right) = D_0.
\]
Concerning the remaining reservoir, we have $\vfinal^{(r)}=V_1^{(r)}=V_0 - C_0$ and $\Dfinal^{(r)}=0.$ Finally, we can approximately determine the amount of $V$ and $D$ coming from infection and lysis processes plus the reservoir. From the process of infection and lysis, from $\CV^{(2)}$-lysis we obtain $\eta (C_0-D_0)$ HV and $\beta \eta (C_0-D_0)$ DVGs. Likewise, from $\CDV$, $\eta D_0/\kappa$ HV are produced and $\eta D_0 (\beta+\delta)/\kappa$ DVGs. Lumping both expressions together, it turns out that the quantity of particles derived from the infection plus lysis process is given by
\begin{eqnarray*}
\vfinal^{(i)} &=& \eta (C_0-D_0) + \frac{\eta}{\kappa} D_0, \\[1.2ex]
\Dfinal^{(i)} &=& \beta \eta (C_0-D_0) + \frac{\eta}{\kappa} (\beta+\delta) D_0.
\end{eqnarray*}
Thus, adding both contributions, one obtains
\begin{eqnarray*}
\vfinal &=& \vfinal^{(i)} + \vfinal^{(r)}= \eta (C_0-D_0) + \frac{\eta}{\kappa} D_0 + (V_0-C_0) \\[1.2ex]
\Dfinal &=& \Dfinal^{(i)} + \Dfinal^{(r)}= \beta \eta (C_0-D_0) + \frac{\eta}{\kappa} (\beta+\delta) D_0,
\end{eqnarray*}
and its sum
\begin{equation*}
\vfinal + \Dfinal = (V_0-C_0) + (1+\beta) \eta (C_0-D_0) + \frac{1+\beta+\delta}{\kappa} \eta D_0 = V_0 + C_0 \left( B - 1\right),
\label{eq:B:moi}
\end{equation*}
where it has been taken into account that $B = (1+\beta) \eta$ and that
\[
\frac{1+\beta+\delta}{\kappa} \eta = (1+\beta) \eta = B.
\]

\end{enumerate}

\subsection{Sensitivity analysis}
\label{se:VE:param}

A crucial piece of information that mathematical models provide, beyond predictions and qualitative interpretations, is to get insights into the relevance of their parameters determining or affecting its dynamics (sensitivity). It is well known in Dynamical Systems Theory that this information can be obtained from the so-called variational equations. In its general form, they supply information about the dependence of any particular solution with respect to its initial conditions. But besides, variational equations with respect to parameters provide knowledge on the sensitivity of such solutions with regard to these parameters. These two complementary sources of information can be both computed using a similar procedure. Let us first introduce the study regarding parameters and show, afterwards, the one corresponding to the initial conditions.

A practical way of getting variational equations in the first case is to differentiate Eqs.~\eqref{eq:1}-\eqref{eq:6} with respect to the chosen parameter, to commute such derivative with the one with respect to the time variable $t$, and to solve the resulting differential equation. As an illustrative example, we detail the case of the variational with respect to the parameter $\beta$. Let us denote by $(C(t;\beta_0), \CV(t;\beta_0),\ldots, D(t;\beta_0))$ the solution of equations~\eqref{eq:1}-\eqref{eq:6} with initial condition $(C_0,\ldots, D_0)$. Let now $(C(t;\beta), \CV(t;\beta),\ldots, D(t;\beta))$ be the solution of the same Cauchy problem but for a parameter $\beta$ close to $\beta_0$. We wonder how close the solution $(C(t;\beta), \ldots, D(t;\beta))$ will evolve (in time) with respect to $(C(t;\beta_0), \ldots, D(t;\beta_0))$. To do it, we compute the variational equation with respect to $\beta$ along $(C(t;\beta_0), \ldots, D(t;\beta_0))$, given by the following system of ODEs: 
\begin{eqnarray*}
\frac{d}{dt} \left( \frac{\partial C}{\partial \beta}\right) &=&
\frac{\partial}{\partial \beta} \dot{C} =
-\infrate (V+D) \cbeta -\infrate C \vbeta -\infrate C \Dbeta, \\
\frac{d}{dt} \left( 
\frac{\partial \CV}{\partial \beta}\right) &=&
\frac{\partial}{\partial \beta} \dot{\CV} =\infrate V \cbeta - (\infrate D+\alpha) \cvbeta + \infrate C \vbeta - \infrate \CV \Dbeta,
\\
\frac{d}{dt} \left( \frac{\partial \CD }{\partial \beta}\right) &=&
\frac{\partial}{\partial \beta} \dot{\CD} =
\infrate D \cbeta -\infrate V \cdbeta - \infrate \CD \vbeta + \infrate C \Dbeta,
\\
\frac{d}{dt} \left( \frac{\partial \CDV }{\partial \beta}\right) &=&
\frac{\partial}{\partial \beta} \dot{C}_{_{DV}} =
\infrate D \cvbeta + \infrate V \cdbeta -\alpha \cdvbeta + \infrate \CD \vbeta + \infrate \CV \Dbeta,
\\
\frac{d}{dt} \left( \frac{\partial V}{\partial \beta}\right) &=&
\frac{\partial}{\partial \beta} \dot{V} = 
\left( \alpha \frac{\partial \eta}{\partial \beta} \left( \CV + \frac{\CDV }{\kappa} \right) - \frac{\alpha \eta}{\kappa^2} \, \frac{\partial \kappa}{\partial \beta} \, \CDV \right) \\
&-& \infrate V \cbeta + \alpha \eta \cvbeta - \infrate V \cdbeta + \frac{\alpha \eta}{\kappa} \cdvbeta - \infrate (C+\CD ) \vbeta,
\\
\frac{d}{dt} \left( \frac{\partial D}{\partial \beta}\right) &=&
\frac{\partial}{\partial \beta} \dot{D} = 
\left( \alpha \left( \eta + \beta \frac{\partial \eta}{\partial \beta} \right) \left( \CV + \frac{\CDV }{\kappa} \right) - \frac{\alpha \beta \eta}{\kappa^2} \CDV \frac{\partial \kappa}{\partial \beta} + \alpha\beta \frac{\partial \eta}{\partial \beta} \, \frac{\CDV }{\kappa} - \frac{\alpha \etaD}{\kappa^2} \, \CDV \frac{\partial \kappa}{\partial \beta} \right) \\
&-& \infrate D \cbeta + (\alpha \beta \eta - \infrate D) \cvbeta + \frac{1}{\kappa} \left( \alpha \beta \eta + \alpha \etaD \right) \cdvbeta - \infrate (C+\CV) \Dbeta,
\end{eqnarray*}
with
\[
\frac{\partial \eta}{\partial \beta}=-\frac{B}{(1+\beta)^2}, \qquad \qquad \frac{\partial \kappa}{\partial \beta} = - \frac{\delta}{(1+\beta)^2},
\]
and where the variables $C,\CV,\ldots,D$ are evaluated at $(t;\beta_0)$ and $\beta=\beta_0$.
The initial conditions of this system are 
\[
\frac{\partial C}{\partial \beta}(0) = \frac{\partial \CV}{\partial \beta}(0)=\frac{\partial \CD}{\partial \beta}(0)=\frac{\partial \CDV}{\partial \beta}(0)=\frac{\partial V}{\partial \beta}(0)=\frac{\partial D}{\partial \beta}(0)=0.
\]
The solutions $\partial C(t)/\partial \beta$, $\partial \CD(t)/\partial \beta$, ... $\partial D(t)/\partial \beta$ of this system provide the (first order) time evolution of the variation of the variables $C, \CV, \ldots, D$ for $\beta$ close to $\beta_0$.
Thus, for instance, in the case of $V(t;\beta)$ we have
\[
V(t;\beta) = V(t;\beta_0) + \frac{\partial V(t;\beta_0)}{\partial \beta} \, (\beta-\beta_0) + \mathcal{O}\left( (\beta-\beta_0)^2\right),
\]
and therefore, for small values of $|\beta-\beta_0|$, we get that
\[
V(t;\beta) \sim V(t;\beta_0) + \frac{\partial V(t;\beta_0)}{\partial \beta} \, (\beta-\beta_0),
\]
where $\frac{\partial V}{\partial \beta}(t;\beta_0)$ has been obtained solving the variational equation above.

\begin{figure}[h!]
\centering
\includegraphics[width=\textwidth]{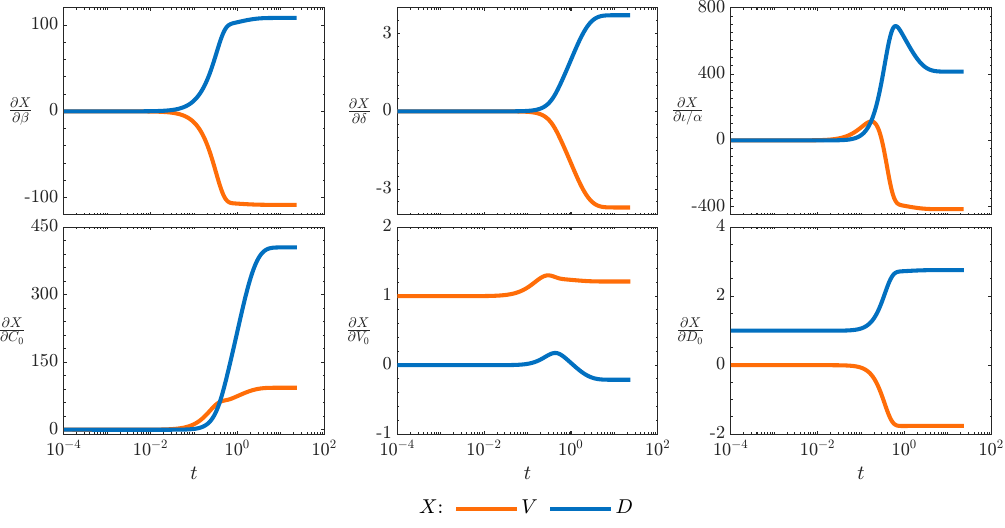}
\caption{Time solutions of the variational equations with respect to parameters $\beta, \delta$, and $\infrate/\alpha$ (top) and initial conditions $C(0),V(0), D(0)$ (bottom) for high multiplicity of infection (MOI) ($m=100$) around the solution plotted in Fig.~\ref{fig:traj_ex} with no degradation and values $qV_0 = 0.75$, $B = 500$, $\beta = 10^{-6}$, $\delta = 10$, and $\infrate/\alpha = 0.1$.} 
\label{fig:VEparam}
\end{figure}

Large values (positive or negative) of $\partial V(t;\beta_0)/\partial \beta$ at a time, say, $t=t^*$ would correspond to significant variation (growth or decay) in the value of $V(t;\beta)$ with respect to $V(t;\beta_0)$ for $\beta \sim \beta_0$. A similar sensitivity analysis can be performed for the rest of the state variables and with respect to the parameters $B, \delta$, and $\infrate/\alpha$. Furthermore, since the experimental data measurements are reduced to values of $C, V$, and $D$, we focus our attention only on the variation of these variables with respect to the selected parameters.
An illustrative example of the solution of these variational equations (with respect to the parameters $B$, $\beta$, $\delta$, and $\infrate/\alpha$) can be found in Fig.~\ref{fig:VEparam}. On one hand, the variationals for $V$ are always in the negative domain, suggesting that increases in three parameters always translate in reductions of HV accumulation. On the other, the variationals for $D$ are always in the positive domain, indicating that increases in the magnitude of any of the three parameters result in more accumulation of DVGs at the cost of the HV. Interestingly, the variational equations with respect to $\delta$ show a small range of values, thus suggesting a weak dependence of the dynamics on the advantage DVGs have on HV in terms of replication. In sharp contrast, the dependence on $\beta$ is strong and even stronger in the case of $\infrate/\alpha$. Together, these observations support the idea that DVGs mostly gain advantage by interfering with the virus (via $\delta$) and, interestingly, reducing the efficiency by which free virus results in virus-producing infected cells (see \eqref{lameua} and discussion therein).

Variational analysis can also be extended to the initial conditions of the system. The derivation of the variational equations of a solution with respect to its initial conditions is defferred to a brief explanation in Section \ref{appendix:se:VE:wrt:ic} in the Appendix. As in here, our analysis only focused on the variables $C, V$, and $D$. Figure~\ref{fig:VEparam} illustrates the effect of $C(0)$, $V(0)$ and $D(0)$ on the result of the process. Here, the situation is more complex than in the case of parameters. The number of available cells at the beginning of the infection has a strong positive effect both in the $V$ and $D$, although the accumulation of DVGs is remarkable more sensitive to $C(0)$ than the accumulation of HVs. In contrast, the effect of the starting number of HVs and of DVGs has minor effects on the outcome of the process. As expected, starting with more $V$ particles slightly benefits the HV and weakly penalizes the DVGs. This has to be understood in a situation in which $m = 100$ and $qV_0 = 0.75$. Likewise, increasing the number of $D$ particles in the inoculum slightly penalizes the HV in the same small magnitude that favors DVGs accumulation.

\subsection{Model fitting to HCoV-OC43 data: parameters' estimation}
\label{subsec:fitting}

To assess the accuracy of the model in reproducing real time-series data and to estimate model parameters, we have fitted the model using artificial intelligence to the data reported in Hillung \emph{et al.} \cite{Hillung2024a}, shown in Fig.~\ref{fig:diagram}b. The aim of such work was to compare the decay of the HCoV-OC43 in cell culture lysates and in fresh media. For this reason the experiment was conducted for several hours after cellular extinction, including degradation processes in our model. A genetic algorithm (GA) has been used to estimate the vector of parameters ($B, \beta, \delta, \alpha, \infrate, \gamma$), and the initial conditions $V_0, D_0$ better fitting the experimental data. The size of the population of vectors of parameters was fixed at $600$, and the number of generations was set to $10^4$. The exploration of the parameter space was refined through successive searches, and the parameter’s range constrained by biological and experimental considerations. Within each generation, the 5\% of the parameters with the lowest $\mathcal{F}$  values (see below) were designated as the elite population, remaining constant for the next generation. The next-best 80\% underwent parameter crossover, while the remaining 15\%, representing the least favorable parameters, experienced random mutations. The parameter space was subdivided into a logarithmic scale to ensure a more balanced exploration, given that the intervals spanned several orders of magnitude. Exceptions were made for $\delta$ and $\gamma$, where linear scale searches were performed due to the nature of their intervals. Concerning the degradation rate, $\gamma$, a narrow interval was chosen around the reported values at \cite{Hillung2024a}. The explored parameter ranges are listed in the Table~\ref{tab:fitting:parameters} below.

\begin{table}[h!]
\caption{Ranges of the parameters of the model explored by the genetic algorithm.}
\centering
\begin{tabular}{rcl} 
\hline\\[-2ex]
 \multicolumn{3}{l}{Parameter} \\ [0.8ex] \hline \\[-1ex]
$B$&$\in$&$[10^1 , 10^4]$ \\ [+0.5ex]
$\beta $&$\in$&$ [10^{-8} , 10^0]$ \\ [+0.5ex]
$\delta $&$\in$&$ [1, 200]$ \\ [+0.5ex]
$\alpha $&$\in$&$ [10^{-3} , 10^{1}] \textrm{ h}^{-1}$ \\ [+0.5ex]
$\infrate $&$\in$&$ [10^{-5} , 10^0] \textrm{ h}^{-1}$ \\ [+0.5ex]
$\gamma  $&$\in$&$ [0.01, 0.10] \textrm{ h}^{-1}$ \\ [+0.5ex]
$V_0 $&$\in$&$ \{[10^5 , 10^7] \textrm{ for } m_V = 1.8 \parallel [10^6 , 10^8] \textrm{ for } m_V = 3.8 \}$ \\ [+0.5ex]
$D_0 $&$\in$&$ \{[1 , 10^7] \textrm{ for } m_V = 1.8 \parallel [1 , 10^8] \textrm{ for } m_V = 3.8 \}$ \\[1.5ex] \hline
\end{tabular}
\label{tab:fitting:parameters}
\end{table}

The cost function used to compute the differences between the experimental and the simulated data is defined as follows:
\begin{equation}
\mathcal{F} = \log \Bigg( \sum_{i=1}^N  \bigg( \log (V_i\cdot \widehat{V}_i^{-1}) \bigg)^2+\sigma \bigg( \log\sum_{j=1}^4\widehat{C}_j(t = 65);k,x^\ast,p \bigg) \Bigg)
\end{equation}
with
\begin{equation}
\sigma(x;k,x^\ast,p)=\frac{p}{1+e^{-k(x-x^\ast)}}.
\end{equation}

Here, $\{t_i,V_i\}_{i = 1}^N$ represent the experimental data  and $\widehat{V}_i$ denotes the  values for the viral population obtained with the mathematical model for a given set of parameters. Additionally, $\sum_j^4 \widehat{C}_j(t=65)$ corresponds to the total number of cells at $t = 65$ hpi, a time point estimated to have a low remaining cell count. The function $\mathcal{F}$ is formulated as the logarithm of a $\chi^2$-function for the logarithm of the experimental data, augmented by a penalty term. The penalty term takes the form of a sigmoidal function with adjustable parameters and is incorporated to account for the approximate time of cell death on the culture, enhancing alignment with the experimental dataset. The parameters of the penalty term have been empirically adjusted to guide the GA in selecting parameter families where the numerical integration of the model significantly reduces the cell count at $t = 65$ hpi (with $k = 5$, $x^\ast = -1$, $p = 10$). External estimates of the cell count over time would enable the removal of this penalty term. To assess the robustness of the final result, it was decided to conduct five identical batches, thereby obtaining five distinct populations of optimized parameters. While performing this procedure, we observed that after $10^4$ generations, the parameters of the elite population were nearly identical. Consequently, we opted to extract the parameter vectors from each batch and compute the mean and standard deviation across batches. 

\begin{figure}
\centering
\includegraphics[]{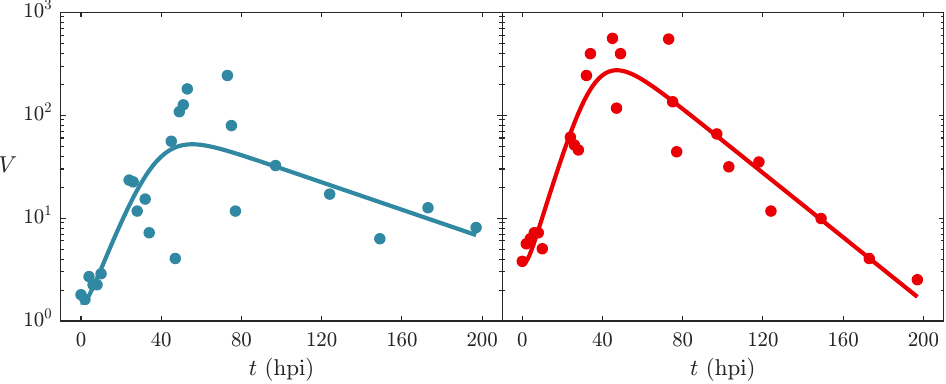}
\caption{Experimental data (circles) and the corresponding fitting of the mathematical model (solid lines) using the parameters values optimised by the genetic algorithm. We show the cases for inoculation conducted at viral multiplicity of infection (MOI) $1.8$ (blue) and $3.8$ (red).} 
\label{fig:fittings}
\end{figure}

For the dataset obtained from inoculation at viral MOI = 1.8, the estimated parameters that best fitted the experimental data were found to be $B = 78$, $\beta = 10^{-8}$, $\delta = 1$, $\alpha = 0.1827 \textrm{ h}^{-1}$, $\infrate = 0.0027 \textrm{ h}^{-1}$ and $\gamma = 0.0154 \textrm{ h}^{-1}$, with corresponding initial conditions $V_0 = 1,516,364$ and $D_0 = 1$. Similarly, for the ones obtained at viral MOI = 3.8, the estimated parameters yielding the optimal reproduction of experimental data were identified as $B = 523$, $\beta = 10^{-8}$, $\delta = 200$, $\alpha = 0.1898 \textrm{ h}^{-1}$, $\infrate = 0.0006 \textrm{ h}^{-1}$ and $\gamma = 0.0361 \textrm{ h}^{-1}$, with corresponding initial conditions $V_0 = 3,676,350$ and $D_0 = 1$. The standard deviation of the optimal parameters in both cases was less than 0.005\% in all cases. The fittings of the mathematical model to the experimental data using the parameters obtained with the GA are displayed in Fig.~\ref{fig:fittings}.

The optimal combination for both datasets, as evident from the resulting optimal combination, shows that $\beta$ and $\delta$ do not find an optimal value within the search range. However, examining the value for $\delta$, it is noteworthy that under such similar experimental conditions, their results are remarkably different. Consequently, we decided to analyze how the function $\mathcal{F}$ behaved while keeping the remaining parameters fixed and varying both $\beta$ and $\delta$ (Appendix Fig.~\ref{fig:app:fittingDetails}). As a result, we found that the $\mathcal{F}$ function hardly changed its value within the parameter search range for $\beta$ and $\delta$. This observation arises because one of the parameter combinations optimizing $\mathcal{F}$ is the absence DVGs (note that $D_0$ does not also find an optimal value within the search range). In the nearly DVGs-free scenario (small $\beta$ and $D_0$), the value of $\delta$ becomes irrelevant.

These results underscore the flexibility of the model and the necessity of acquiring additional data for the fitting process for a rigorous interpretation of the parameters. Measures that could prove helpful would include, for instance, counting total cell numbers over time using life imaging techniques and teasing apart infected from non-infected cells by using viruses tagged with fluorescent proteins and living from dead cells by using dead-specific dyes such as trypan blue or propidium iodide. Lastly, to precisely estimate the amount of DVGs, it would be essential to extract total RNA from the cultures at various time points, subjecting them to RNA-seq analysis followed by bioinformatic analyses using tools such as DVGfinder \cite{OlmoUceda2022} to identify different DVGs and estimators as described in \cite{MunozSanchez2024} to quantify their abundance.

\section{Discussion}

Defective viral genomes (DVGs) have been often viewed as replication byproducts, stemming from cell culture passaging conditions. However, a renewed interest in this fraction of the total viral population suggests that their existence might have more functional or biological significance than previously thought~\cite{Vignuzzi2019}. Given the notably high recombination rates seen in some viruses, including betacoronaviruses, it is plausible that DVGs could influence the adaptive evolution of the viable virus population~\cite{Campos2022}. The question of whether DVGs confer or not a selective advantage to the viral population and if the conditions of high multiplicity of infection (MOI) or localized co-infection hold biological relevance beyond cell culture remains to be explored.

Serial passage experiments serve for many different applications in Virology, not only studying \emph{in vitro} evolution \cite{Elena2007}. The two most common applications are (\emph{i}) the maintenance and amplification of viral stocks in the laboratory and (\emph{ii}) the classic way of producing attenuated life vaccines by adapting the target virus to cells from a host as different as possible from the one to be immunised. DVGs have been found in live-attenuated vaccines for polio, measles and influenza viruses (reviewed in \cite{Vignuzzi2019}). However, their impact on the development of protective immunity and vaccine efficacy has not been formally evaluated. Given their potential to interfere with and stimulate the immune system, there is speculation that DVGs could improve vaccine efficiency while ensuring virus safety by limiting its replication and spread. If this hypothesis holds true, it becomes crucial to carefully control the amount of DVGs in vaccine preparations to prevent complete interference and a significant reduction in the virus' effectiveness. Our results suggest that the number of DVGs generated from helper virus (HV) replication ($\beta$) and the replicative advantage of DVGs ($\delta$) are the two more relevant parameters to be manipulated in order to optimize the ratio HV:DVGs in terms of vaccine efficiency.

Multitude of mathematical models for HV-DVGs have been investigated at different scales \cite{Gao2009, Bangham1990, Szathmary1992, Szathmary1993, Frank2000, Sardanyes2010, StaufferThompson2010, Zwart2013, Chao2017, Chaturvedi2021, Liang2023}. Interestingly, these models have reproduced experimentally observed dynamics such as deterministic chaos \cite{Zwart2013} and self-curing \cite{Kirkwood1994, StaufferThompson2010}, which involves the simultaneous extinction of the HV and defective interfering particles (DIPs), result that was observed \emph{in vitro} by Jacobson \textit{et al.} \cite{Jacobson1979} and Stauffer Thompson and Yin \cite{StaufferThompson2010}. Our model also identifies combinations of parameters in which self-curing arises: large $\delta$ values would in turn allow for the relationship $1 + \beta + \delta > B$ to be fulfilled, shifting from the scenario of a single attracting plane $\PiVD$ to a more complex scenario with three attracting planes and tristability, including the planes $\PiCCD$ and $\PiCDD$ (being these planes degenerate NHIMs) representing virus-free self-curing solutions.

Model predictions are as valid and realistic as the assumptions from which they build up. Our first strong assumption, already discussed at large, was that at the time scales of the evolution experiments, virus production over-weights virus degradation. This assumption is well supported by experimental data (see~\cite{Hillung2024b} and references therein as well as \cite{Martinez2011}). Our model was designed to shed light on the dynamics of betacoronavirus populations during serial passages and thus could not be applied to a much more complex \emph{in vivo} situation. For example, we are collapsing all possible DVGs into a single category, while in reality, viral populations contain a large fraction of very diverse DVGs that are in a dynamic equilibrium, with some appearing and disappearing at every transmission event and others persisting for long periods of time \cite{Jaworski2017, Rangel2023, Zhou2023, Hillung2024a}. Our approach aligns with the current mathematical models mentioned in the previous paragraph that typically focus on a single dominant DVG. Further developments to incorporate the potential cooperation and competition among multiple types of DVGs are needed.

Our model is similar to those proposed by Frank \cite{Frank2000} and by Liang \emph{et al.} \cite{Liang2023} in the sense that they incorporate $V$, $D$, $C$, $\CD$, $\CDV$, and $\CV$ as state variables. However, there are two major differences. Firstly, their models consider a spatial distribution of cells and thus use reaction-diffusion systems. Secondly, they incorporate an additional category of cells, namely $\CV^*$, considered as cells that were infected by $V$ some time ago, becoming $\CV$, and are already producing HV but cannot be superinfected by DVGs. In this case, the fate of $\CV^*$ is to die out after some time producing only $V$. This makes a notable difference with our model, as we opted for simplicity by suppressing this category as the modulation of the quotient $\infrate/\alpha$ embeds the scenario in which $\CV$ lysates previous to superinfection with $D$ (Appendix Fig.~\ref{fig:app:ia_figure}). Although we recognize that superinfection exclusion \cite{Hunter2022} might be a relevant process for many viruses, it is not a universal one, and no evidence of such process have been reported for coronaviruses.

The model presented here provides two new results in theoretical Virology: (\emph{i}) HV-DVGs replication and infection dynamics governed by quasineutral planes (degenerate 2D NHIMs formed by equilibrium points in this case); (\emph{ii}) scenario of tristability formed by these degenerate NHIMs. Concerning point (\emph{i}) and, as far as we know, quasineutral dynamics in RNA virus dynamics have been up to now found in degenerate one-dimensional objects~\cite{Sardanyes2018}. Hence, our findings extend this result on degenerate manifolds to planes. Concerning point (\emph{ii}), these mentioned one-dimensional manifolds are usually found in scenarios of monostability and bistability. We here provide an example where three different quasineutral planes can be achieved depending on initial conditions. The existence of these degenerate planes is subject to non spontaneous degradation of viral particles and cells ($\gamma = 0$), scenario that may seem unlikely in a real experiment. However, experimental estimates have indicated that degradation remains very low and that the dynamics within the time scale used in the experiments is dominated by virus replication over degradation~\cite{Hillung2024b}. Hence, we conjecture that the long time delays arising with parameter values close to the ones at which these degenerate NHIMs are found~\cite{Fontich2022} could be observed in real experiments. Moreover, the passage experiments showed extremely large fluctuations between passages (some of them of about 1-2 orders of magnitude), as we show in Fig.~\ref{fig:diagram}b (see~\cite{Hillung2024a} for further details). These extremely large fluctuations could be due to the combined effect of stochastic sampling between passages and the dynamics tied to degenerate NHIMs. Under this scenario, different initial conditions (stochastically varying at each passage) may give place to different stationary values of HV and DVGs due to influence of NHIMs or to their remnants for those cases with some (and small) contribution of degradation. 

In conclusion, the mathematical model studied in this manuscript fits well the dynamics of virus accumulation within single-passages \emph{i.e.}, within- and between-cell dynamics in a cell culture, in cell cultures of betacoronaviruses. Our parameter sensitivity analysis also suggests that the most relevant parameters to explain the observed dynamical patterns are the production of DVGs from the HV ($\beta$) and the infection to lysate rates ($\infrate/\alpha$). The combination of experiments and biologically inspired modelling provides a powerful tool to identify parameters to be optimized for future developments of DVGs as therapeutic interfering particles or as adjuvants in life attenuated vaccines.

\section*{Acknowledgements}

We thank José A. Oteo for valuable suggestions and critical reading of the mannuscript. JCM has been funded by grant ACIF/2021/296 (Generalitat Valenciana). MJO was funded by contract FPU19/05246 by MCIU/AEI/ 10.13039/501100011033 and ‘‘ESF invests in your future’’. JTL has been funded by the projects PGC2018-098676-B-100 and PID2021-122954NB-I00 funded by MCIU/AEI/10.13039/501100011033/ and ‘‘ERDF a way of making Europe’’, and by the grant ‘‘Ayudas para la Recualificaci\'on del Sistema Universitario Espa\~nol 2021-2023’’. JTL also thanks the Laboratorio Subterr\'aneo de Canfranc, the I$^2$SysBio and the Institut de Math\'ematiques de Jussieu-Paris Rive Gauche (Sorbonne Universit\'e) for their hospitality as hosting institutions of this grant. 
We also thank the AEI, through the Mar\'ia de Maeztu Program for Units of Excellence in R\&D (CEX2020-001084-M) and CERCA Programme/Generalitat de Catalunya for institutional support. JS has been also supported by the Ram\'on y Cajal grant RYC-2017-22243 funded by MCIU/AEI/10.13039/501100011033 and ‘‘ESF invests in your future’’.  SFE was supported by CSIC PTI Salud Global grant 202020E153 and by grants SGL2021-03-009 and SGL2021-03-052 from European Union Next Generation EU/PRTR through the CSIC Global Health Platform established by EU Council Regulation 2020/2094. Many computations were performed on the HPC cluster Garnatxa at I$^2$SysBio (CSIC-UV).  JCMS, JTL, MJOU, and SFE acknowledge the support of the Santa Fe Institute, where part of this research was developed.


\newpage
\begin{center}
{\Large{\textsc{Appendix}}    }
\end{center}

\appendix
\setcounter{equation}{0}
\renewcommand{\theequation}{\thesection\arabic{equation}}
\setcounter{figure}{0}
\renewcommand{\thefigure}{\thesection\arabic{figure}}
\setcounter{table}{0}
\renewcommand{\thetable}{\thesection\arabic{table}}

\section{Time-series examples for the $\PiVD$, $\PiCDD$ and $\PiCCD$}
\label{appendix:timeSeries_examples}
\setcounter{figure}{0}

\begin{figure}[h!]
\centering
\includegraphics[width=\textwidth]{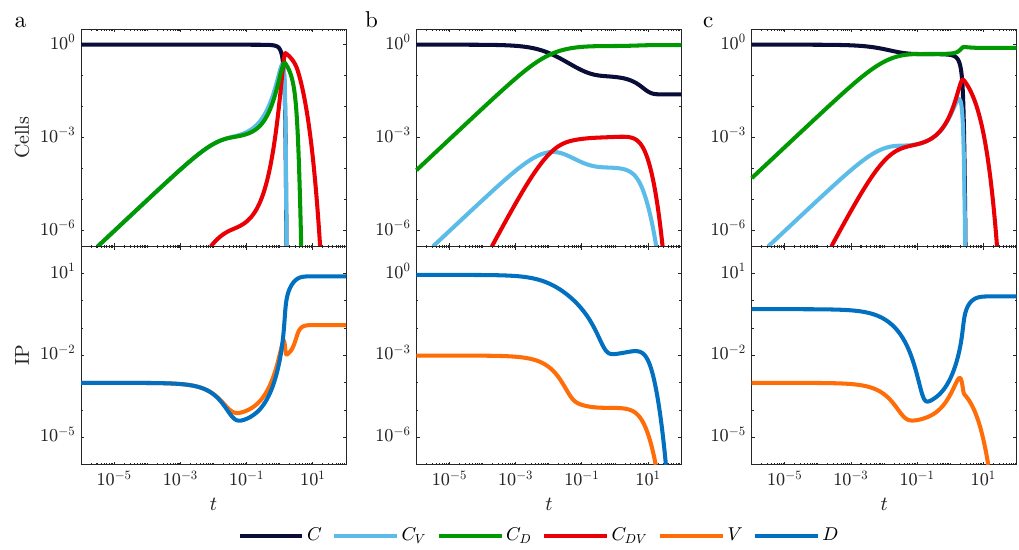}
\caption{Time series to extinctions for cells (end state at \textcolor{PiVD_color}{$\PiVD$}) and for viral infective particles (IP) (end state at \textcolor{PiCdD_color}{$\PiCDD$}) and double extinction (end state \textcolor{PiCCd_color}{$\PiCCD$}) Parameters: $B = 10$, $\beta = 0.5$ and $\delta = 20$ and $\infrate/\alpha = 100$. Initial conditions: (a) $V_0 = 10^{-3}$ and $D_0 = 10^{-3}$; (b) $V_0 = 10^{-3}$ and $D_0 = 0.9$; and (c) $V_0 = 10^{-3}$ and $D_0 = 0.5$. $C_0 = 1$.}
\label{fig:app:comp:triestability_timeSeries}
\end{figure} 

\newpage
\section{First order variational equation with respect to initial conditions}
\label{appendix:se:VE:wrt:ic}
In this section we present a well-known derivation of the so-called first variational equation with respect to initial conditions. To this end, let us consider the following Cauchy problem:
\begin{equation}
\dot{x}=f(x), \qquad x(0)=x_0, \quad \qquad x\in \R^n, \quad \dot{\phantom{x}}=\frac{d}{dt},
\label{eq:VE:init}
\end{equation}
which, without loss of generality, we can assume autonomous and satisfying all the required conditions of smoothness and derivability, ensuring the existence and uniqueness of solutions. We denote by $x(t;x_0)$ its unique solution. We are interested in the computation of the first-order (in $\eps$) variation of the solution $x(t;x_0+\eps)$ of the Cauchy problem 
\begin{equation}
\dot{x}=f(x), \qquad \qquad x(0)=x_0+\eps, 
\label{eq:VE:pert}
\end{equation}
where, abusing of notation, $x_0+\eps= x_0 + \eps \mathrm{Id}$, $\mathrm{Id}$ being the $n$-dimensional identity matrix. If $\eps$ is small enough, we can write, using Taylor,
\begin{equation}
x(t;x_0+\eps) = x(t;x_0) + \eps \frac{\partial x(t;x_0)}{\partial x_0} + \mathcal{O}(\eps^2).
\label{eq:VE:taylor}
\end{equation}
Differentiating this expression with respect to $t$ we get
\begin{equation}
\dot{x}(t;x_0+\eps) = \dot{x}(t;x_0) + \eps \frac{d}{dt} \left( \frac{\partial x(t;x_0)}{\partial x_0} \right) + \mathcal{O}(\eps^2).
\label{eq:VE:1}
\end{equation}
On the other side, using that $x(t;x_0+\eps)$ is a solution of~\eqref{eq:VE:pert} and using Taylor again, we have that
\begin{eqnarray*}
\dot{x}(t;x_0+\eps) &=& f(x(t;x_0+\eps)) = f(x(t;x_0)) + 
\eps \frac{\partial f(x(t;x_0))}{\partial x_0} +  \mathcal{O}(\eps^2) \\
&=& f(x(t;x_0)) + 
\eps Df(x(t;x_0)) \, \frac{\partial x(t;x_0)}{\partial x_0} + \mathcal{O}(\eps^2),
\end{eqnarray*}
where $Df$ denotes the differential matrix of $f$ with respect to $x=(x_1,x_2,\ldots,x_n)$. Equating terms of order $\eps^1$ in the latter expression with those in~\eqref{eq:VE:1} it turns out that $\partial x(t;x_0)/\partial x_0$ satisfies the ODE
\begin{equation}
\frac{d}{dt} \left( \frac{\partial x(t;x_0)}{\partial x_0} \right) =
Df(x(t;x_0)) \, \frac{\partial x(t;x_0)}{\partial x_0}.
\label{VE}
\end{equation}
Regarding its initial condition, we know that (substituting $t=0$ in expression~\eqref{eq:VE:taylor})
\[
x(0;x_0+\eps)=x(0;x_0) + \eps \frac{\partial x(0;x_0)}{\partial x_0} + \mathcal{O}(\eps^2) \Leftrightarrow
x_0+\eps=x_0 + \eps \frac{\partial x(0;x_0)}{\partial x_0} + \mathcal{O}(\eps^2),
\]
from where, equating again powers in $\eps$, it follows that
\begin{equation}
\frac{\partial x(0;x_0)}{\partial x_0} = \mathrm{Id}, 
\label{VE:ic}
\end{equation}
where $\mathrm{Id}$ is the identity matrix. Equations~\eqref{VE} and~\eqref{VE:ic} are usually called first variational equations around a solution $x(t;x_0)$. They provide information about the local dynamics (tangential and normal) around a solution $x(t;x_0)$ of~\eqref{eq:VE:init}.   
This first variational equation is linear and homogeneous. Recurrently, one can compute the variational equations of higher order (in $\eps$). All of them are also linear but non-homogeneous.

\newpage
\section{Additional fittings results}
\label{appendix:fittingDetails}
\setcounter{figure}{0}

\begin{figure}[!ht]
\centering
\includegraphics[width = \textwidth]{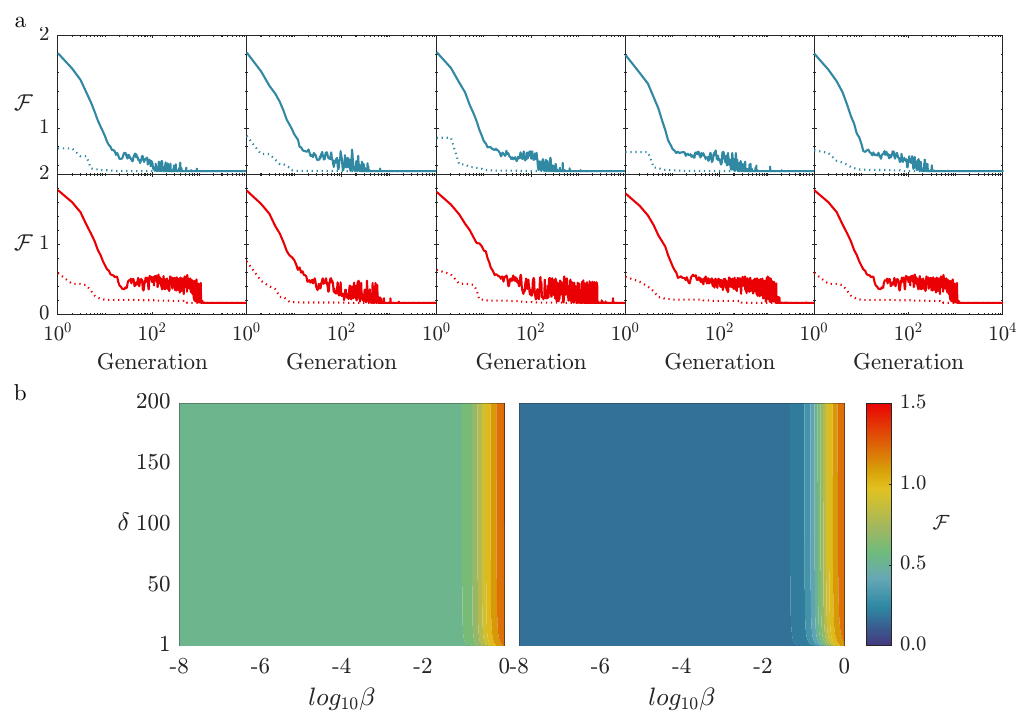}
\caption{(a) Mean (solid line) and minimum (doted line) value of $\mathcal{F}$ at each generation for data from inoculation at viral multiplicity of infection (MOI) $1.8$ (blue) and $3.8$ (red). (b) Value of $\mathcal{F}$ for different ($\beta$, $\delta$) combinations being the rest of the parameters selected as optimum at the genetic algorithm. The same study has been performed with $D_0 = 0$ for both datasets obtaining the same results.} 
\label{fig:app:fittingDetails}
\end{figure}

\newpage
\section{Model solutions for some $\infrate/\alpha$ values}

\begin{figure}[h!]
\centering
\includegraphics[width=\textwidth]{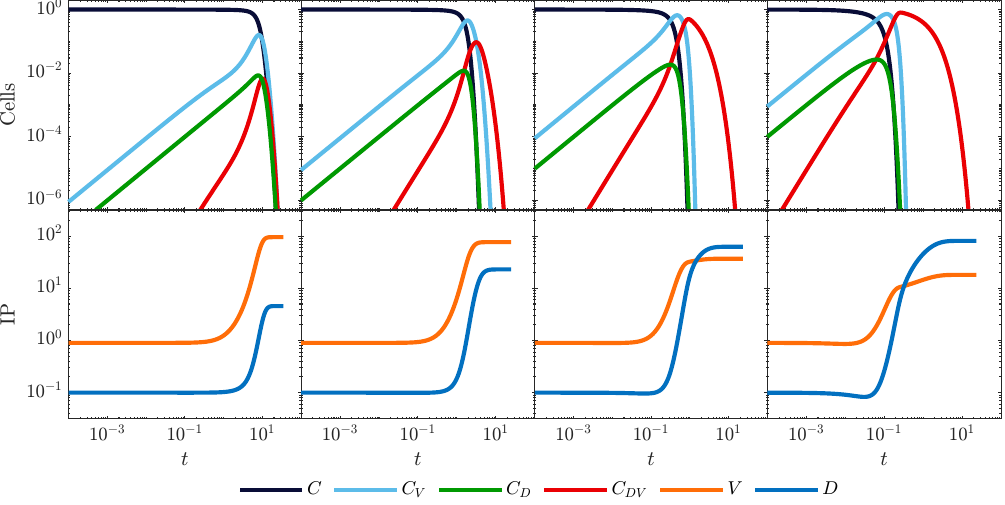}
\caption{Some numerical integrations of the systems for $B = 100$, $\beta = 0.01$ and $\delta = 10$. From left to right the $\infrate/\alpha$ quotient is 0.01, 0.1, 1 and 10. This figure aims to illustrate how the result depends on these parameters, both in terms of cellular type abundance and the helper virus (HV) - defective viral genomes (DVGs) final ratio (IP: infective particles).}
\label{fig:app:ia_figure}
\end{figure}

\end{document}